\documentclass[11pt]{amsart}%
\usepackage{amsfonts}
\usepackage{amsmath}
\usepackage{amssymb}
\usepackage{graphicx}%
\newtheorem{theorem}{Theorem}

\newtheorem{corollary}[theorem]{Corollary}

\newtheorem{definition}[theorem]{Definition}
\newtheorem{example}[theorem]{Example}

\newtheorem{lemma}[theorem]{Lemma}

\newtheorem{proposition}[theorem]{Proposition}

\newenvironment{notation}{\noindent {\bf Notation }}{\\}

\newcommand{\E}{{\mathcal{E}}}
\begin{document}
\title[Mixed Tsirelson Spaces]{More Mixed Tsirelson Spaces That Are Not Isomorphic To Their Modified Versions}
\author{Denny H. Leung}
\address{Department of Mathematics, National University of Singapore, 2 Science Drive
2, Singapore 117543.}
\email{matlhh@nus.edu.sg}
\author{Wee-Kee Tang}
\address{Mathematics and Mathematics Education, National Institute of Education \\
Nanyang Technological University, 1 Nanyang Walk, Singapore 637616.}
\email{wktang@nie.edu.sg}

\begin{abstract}
The class of mixed Tsirelson spaces is an important source of examples in the
recent development of the structure theory of Banach spaces. The related class
of modified mixed Tsirelson spaces has also been well studied. In the present
paper, we investigate the problem of comparing isomorphically the mixed
Tsirelson space $T[({\mathcal{S}}_{n},\theta_{n})_{n=1}^{\infty}]$ and its
modified version $T_{M}[({\mathcal{S}}_{n},\theta_{n})_{n=1}^{\infty}]$. It is
shown that these spaces are not isomorphic for a large class of parameters
$(\theta_{n})$.

\end{abstract}
\maketitle

\section{Introduction}

In 1974, Tsirelson \cite{T} settled a fundamental problem in the structure
theory of Banach spaces when he gave a surprisingly simple construction of a
Banach space that does not contain any isomorphic copy of $c_{0}$ or $\ell
^{p}$, $1\leq p<\infty$. Figiel and Johnson \cite{FJ} provided an analytic
description, based on iteration, of the norm of the dual of Tsirelson's
original space. Subsequently, other examples of spaces were constructed with
norms described iteratively, notable among them were Tzafriri's spaces
\cite{Tz} and Schlumprecht's space\cite{S}. Gowers' and Maurey's solution to
the unconditional basic sequence problem\ \cite{GM}\ is a variation based on
the same theme. It has emerged in recent years that, far from being isolated
examples, Tsirelson's space and its variants from an important class of Banach
spaces. Argyros and Deliyanni \cite{AD} were the first to provide a general
framework for such spaces by defining the class of mixed Tsirelson spaces.
Among the earliest variants of Tsirelson's space was its modified version
introduced by Johnson \cite{J}. Casazza and Odell \cite{CO} showed that
Tsirelson's space is isomorphic to its modified version. This isomorphism was
exploited to study the structure of the space. The modification can be
extended directly to the class of mixed Tsirelson spaces, forming the class of
modified mixed Tsirelson spaces. It is thus of natural interest to determine
if a mixed Tsirelson space is isomorphic to its modified version. This
question has been considered by various authors, e.g., \cite{ADKM, M}, who
provided answers in what may be considered \textquotedblleft
extremal\textquotedblright\ cases. In the present paper, we show that for a
large class of parameters, a mixed Tsirelson space and its modified version
are not isomorphic.

We shall be concerned exclusively with mixed Tsirelson spaces of the form
$T[({\mathcal{S}}_{n},\theta_{n})^{\infty}_{n=1}]$ or $T[({\mathcal{S}}%
_{n_{i}},\theta_{i})^{k}_{i=1}]$ and their modified versions. We now recall
the definitions of these spaces and the various notions involved. Denote by
${\mathbb{N}}$ the set of natural numbers. For any infinite subset $M$ of
${\mathbb{N}}$, let $[M]$ and $[M]^{<\infty}$ be the set of all infinite and
finite subsets of $M$ respectively. These are subspaces of the power set of
${\mathbb{N}}$, which is identified with $2^{{\mathbb{N}}}$ and endowed with
the topology of pointwise convergence. If $I$ and $J$ are nonempty finite
subsets of ${\mathbb{N}}$, we write $I<J$ to mean $\max I<\min J$. We also
allow that $\emptyset<I$ and $I<\emptyset$. For a singleton $\{n\}$, $\{n\}<J$
is abbreviated to $n<J$. The general Schreier families ${\mathcal{S}}_{\alpha
}$, $\alpha< \omega_{1}$, were introduced by Alspach and Argyros \cite{AA}. We
shall restrict ourselves to finite parameters. Let ${\mathcal{S}}_{0}$ consist
of all singleton subsets of ${\mathbb{N}}$ together with the empty set.
Inductively, if $n \in{\mathbb{N}}$, let ${\mathcal{S}}_{n}$ consist of all
sets of the form $\cup^{k}_{i=1}G_{i}$, where $G_{i} \in{\mathcal{S}}_{n-1}$,
$G_{1} < \dots< G_{k}$ and $k \leq\min G_{1}$. The Schreier families are
\textbf{hereditary}: $G \in{\mathcal{S}}_{n}$ whenever $G \subseteq F$ and $F
\in{\mathcal{S}}_{n}$; \textbf{spreading}: for all strictly increasing
sequences $(m_{i})_{i=1}^{k}$ and $(n_{i})_{i=1}^{k}$, $(n_{i})_{i=1}^{k}%
\in{\mathcal{S}_{n}}$ if $(m_{i})_{i=1}^{k}\in{\mathcal{S}_{n}}$ and
$m_{i}\leq n_{i}$ for all $i$; and compact as subspaces of $[{\mathbb{N}%
}]^{<\infty}$. A sequence in $(E_{i})^{k}_{i=1}$ in $[{\mathbb{N}}]^{<\infty}$
is said to be ${\mathcal{S}}_{n}$\textbf{-admissible} if $E_{1} < \dots<
E_{k}$ and $\{\min E_{i}\}^{k}_{i=1} \in{\mathcal{S}}_{n}$. It is
${\mathcal{S}}_{n}$\textbf{-allowable} if the $E_{i}$'s are pairwise disjoint
and $\{\min E_{i}\}^{k}_{i=1} \in{\mathcal{S}}_{n}$.

Denote by $c_{00}$ the space of all finitely supported real sequences, whose
unit vector basis will be denoted by $(e_{k})$. For a finite subset $E$ of
${\mathbb{N}}$ and $x\in c_{00}$, let $Ex$ be the coordinatewise product of
$x$ with the characteristic function of $E$. The $\sup$ norm and the $\ell
^{1}$-norm on $c_{00}$ are denoted by $\Vert\cdot\Vert_{c_{0}}$ and
$\Vert\cdot\Vert_{\ell^{1}}$ respectively. Given a null sequence $(\theta
_{n})_{n=1}^{\infty}$ in $(0,1)$, define sequences of norms $\Vert\cdot
\Vert_{m}$ and $|||\cdot|||_{m}$ on $c_{00}$ as follows. Let $\Vert x\Vert
_{0}=|||x|||_{0}=\Vert x\Vert_{c_{0}}$ and
\begin{equation}
\Vert x\Vert_{m+1}=\max\{\Vert x\Vert_{m},\sup_{n}\theta_{n}\sup\sum_{i=1}%
^{r}\Vert E_{i}x\Vert_{m}\}, \label{eq 0}%
\end{equation}
where the last sup is taken over all ${\mathcal{S}}_{n}$-admissible sequences
$(E_{i})_{i=1}^{r}$. The norm $|||x|||_{m}$ is defined as in ($\ref{eq 0}$)
except that the last sup is taken over all $\mathcal{S}_{n}$-allowable
sequences $\left(  E_{i}\right)  _{i=1}^{r}.$ Since these norms are all
dominated by the $\ell^{1}$-norm, $\Vert x\Vert=\lim_{m}\Vert x\Vert_{m}$ and
$|||x|||=\lim_{m}|||x|||_{m}$ exist and are norms on $c_{00}$. The
\textbf{mixed Tsirelson space} $T[({\mathcal{S}}_{n},\theta_{n})_{n=1}%
^{\infty}]$ and the \textbf{modified mixed Tsirelson space} $T_{M}%
[({\mathcal{S}}_{n},\theta_{n})_{n=1}^{\infty}]$ are the completions of
$c_{00}$ with respect to the norms $\left\Vert \cdot\right\Vert $ and
$|||\cdot|||$ respectively. From equation (\ref{eq 0}) we can deduce that
these norms satisfy the implicit equations%
\begin{equation}
\Vert x\Vert=\max\{\Vert x\Vert_{c_{0}},\sup_{n}\theta_{n}\sup\sum_{i=1}%
^{r}\Vert E_{i}x\Vert\} \label{equation 0.1}%
\end{equation}
and
\begin{equation}
|||x|||=\max\{\Vert x\Vert_{c_{0}},\sup_{n}\theta_{n}\sup\sum_{i=1}%
^{r}|||E_{i}x|||\}, \label{equation 0.2}%
\end{equation}
where the innermost suprema are taken over all ${\mathcal{S}}_{n}$-admissible,
respectively, ${\mathcal{S}}_{n}$-allowable sequences $(E_{i})_{i=1}^{r}$. The
mixed Tsirelson space $T[(\mathcal{S}_{n_{i}},\theta_{i})_{i=1}^{k}]$ and
modified mixed Tsirelson space $T_{M}[(\mathcal{S}_{n_{i}},\theta_{i}%
)_{i=1}^{k}]$ are defined similarly.

For standard Banach space terminology and notation, we refer to \cite{LinT}.
Two Banach spaces $X$ and $Y$ are said to be \textbf{isomorphic} if they are
linearly homeomorphic. A linear homeomorphism from $X$ into $Y$ is called an
\textbf{embedding}. We say that $X$ \textbf{embeds into} $Y$ if such an
embedding exists. $X$ and $Y$ are \textbf{totally incomparable} if no infinite
dimensional subspace of one embeds into the other. A sequence $(x_{n})$ in $X$
is said to \textbf{dominate} a sequence $(y_{n})$ in $Y$ if there is a finite
constant $K$ such that $\|\sum a_{n}y_{n}\| \leq K\|\sum a_{n}x_{n}\|$ for all
$(a_{n}) \in c_{00}$. Two sequences are \textbf{equivalent} if they dominate
each other.

\section{Brief Survey of Known Results}

The aim of the present paper is to compare isomorphically the spaces
$T[({\mathcal{S}}_{n},\theta_{n})_{n=1}^{\infty}]$ and $T_{M}[({\mathcal{S}%
}_{n},\theta_{n})_{n=1}^{\infty}]$ (and also the spaces $T[(\mathcal{S}%
_{n_{i}},\theta_{i})_{i=1}^{k}]$ and $T_{M}[(\mathcal{S}_{n_{i}},\theta
_{i})_{i=1}^{k}]$). Let us recall some known results in this direction.
Casazza and Odell \cite{CO} showed that the Tsirelson space $T[\mathcal{S}%
_{1},\theta]$ is isomorphic to the modified Tsirelson space $T_{M}%
[\mathcal{S}_{1},\theta],$ with no specific isomorphism constant given in
their proof. In \cite{B}, Bellenot proved that they are $\theta^{-1}%
$-isomorphic. Recently, Manoussakis \cite{M} showed that the spaces
$T[\mathcal{S}_{n},\theta]$ and $T_{M}[\mathcal{S}_{n},\theta]$ are
$3$-isomorphic for all $n\in\mathbb{N}$ and all $\theta\in(0,1)$. He also
stated without proof in \cite[Section 4]{M2} that $T[(\mathcal{S}_{n_{i}%
},\theta_{i})_{i=1}^{k}]$ is isomorphic to $T_{M}[(\mathcal{S}_{n_{i}}%
,\theta_{i})_{i=1}^{k}]$. A proof of a nominally more general fact will be
given below.

When considering the spaces $T[({\mathcal{S}}_{n},\theta_{n})^{\infty}_{n=1}]$
and $T_{M}[({\mathcal{S}}_{n},\theta_{n})^{\infty}_{n=1}]$, we may assume
without loss of generality that $(\theta_{n})$ is nonincreasing and that
$\theta_{m+n} \geq\theta_{m}\theta_{n}$. Such sequences are said to be
\textbf{regular}. It is known that \cite{OTW} $\lim\theta_{n}^{1/n} =
\sup\theta^{1/n}_{n}$ for a regular sequence $(\theta_{n})$. Argyros et. al.
showed that if $(\theta_{n})$ is regular and $\lim_{n}\theta_{n}^{1/n}=1,$
then $T[(\mathcal{S}_{n},\theta_{n})_{n=1}^{\infty}]$ contains copies of
$\ell^{\infty}(n)$'s uniformly and hereditarily \cite[Theorem 1.6]{ADKM}. As a
result, they were able to conclude that $T[(\mathcal{S}_{n},\theta_{n}%
)_{n=1}^{\infty}]$ and $T_{M}[(\mathcal{S}_{n},\theta_{n})_{n=1}^{\infty}]$
are totally incomparable.

In \cite{LT2}, the authors introduced the condition
\[
(\dag) \quad\quad\lim_{m}\limsup_{n}\frac{\theta_{m+n}}{\theta_{n}}=1.
\]
Condition $(\dag)$ is weaker than the condition $\lim_{n}\theta_{n}^{1/n}=1$.
Indeed, if
\[
\lim_{m}\limsup_{n}\frac{\theta_{m+n}}{\theta_{n}}<1,
\]
then there exist $\delta<1, m\in\mathbb{N}$ and $N\in\mathbb{N}$ such that
$\frac{\theta_{n+m}}{\theta_{n}}<\delta$ for all $n\geq N.$ In particular, for
all $k\in\mathbb{N}$, $\theta_{km+N}<\delta^{k}\theta_{N}.$ Thus
$\theta_{km+N}^{\frac{1}{km+N}}<\delta^{\frac{k}{km+N}}\theta_{N}^{^{\frac
{k}{km+N}}}.$ Taking $k\rightarrow\infty$, we have $\lim_{n}\theta_{n}%
^{1/n}\leq\delta^{1/m}<1.$ It can be shown that the converse is false, even
for regular sequences.

If $(\theta_{n})$ satisfies (\dag), it follows from \cite[Proposition 9]{LT1}
that there exists $\varepsilon> 0$ such that for all $V\in[\mathbb{N}]$ and
all $k\in\mathbb{N}$, there exists a sequence of pairwise disjoint vectors
$(y_{j}) _{j=1}^{k}$ $\subseteq\operatorname*{span}\left\{  e_{k}:k\in
V\right\}  $ such that $\|\sum_{j=1}^{k}y_{j}\| \leq2 + 1/\varepsilon$ and
$\|y_{j}\| \geq1$ for all $j.$ In other words, $\ell^{\infty}\left(  n\right)
$'s uniformly disjointly embeds into the subspace of $T\left[  \left(
\mathcal{S}_{n},\theta_{n}\right)  _{n=1}^{\infty}\right]  $ generated by
$\left(  e_{k}\right)  _{k\in V}$. In particular, the norms $\|\cdot\|$ and
$|||\cdot|||$ are not equivalent on $\operatorname*{span}\left\{  e_{k}:k\in
V\right\}  $. This together with the proposition below imply that $T\left[
\left(  \mathcal{S}_{n},\theta_{n}\right)  _{n=1}^{\infty}\right]  $ is not
isomorphic to $T_{M}\left[  \left(  \mathcal{S}_{n},\theta_{n}\right)
_{n=1}^{\infty}\right]  .$

\begin{proposition}
\label{Nonisomorphic}If $T[(\mathcal{S}_{n},\theta_{n})_{n=1}^{\infty}]$
embeds into $T_{M}[(\mathcal{S}_{n},\theta_{n})_{n=1}^{\infty}]$, then there
exists $V\in[\mathbb{N}]$ such that $|||\cdot||| $ is equivalent to
$\|\cdot\|$ on the subspace $\operatorname*{span}\{e_{k}: k \in V\}.$
\end{proposition}

\begin{proof}
Let $J:T\left[  \left(  \mathcal{S}_{n},\theta_{n}\right)  _{n=1}^{\infty
}\right]  \rightarrow T_{M}\left[  \left(  \mathcal{S}_{n},\theta_{n}\right)
_{n=1}^{\infty}\right]  $ be an embedding. Then $\left(  Je_{k}\right)  $ is a
weakly null sequence. By the Bessaga-Pe\l czynski Selection Principle (see
e.g. \cite[Proposition 1.a.12]{LinT}), there is a subsequence $(Je_{k_{j}})$
of $\left(  Je_{k}\right)  $ such that $(Je_{k_{j}})$ is equivalent to a
seminormalized block sequence $\left(  u_{j}\right)  $ in $T_{M}\left[
\left(  \mathcal{S}_{n},\theta_{n}\right)  _{n=1}^{\infty}\right]  .$ Let
$m_{j}=\min\operatorname*{supp}u_{j}.$ By taking a subsequence if necessary,
we may assume that $k_{j}\leq m_{j}$ for all $j.$ Choose $j_{1}<j_{2}<\cdots$
such that%
\[
k_{j_{1}}\leq m_{j_{1}}\leq k_{j_{2}}\leq m_{j_{2}}\leq\cdots
\]
Note that in $T_{M}\left[  \left(  \mathcal{S}_{n},\theta_{n}\right)
_{n=1}^{\infty}\right]  ,$ $(e_{m_{j_{i}}})$ is dominated by $\left(
u_{j_{i}}\right)  ,$ which is equivalent to $(Je_{k_{j_{i}}}).$ Hence there
exists a finite constant $\lambda$ such that for all $(a_{i})\in c_{00},$
\begin{align*}
|||\sum a_{i}e_{m_{j_{i}}}|||  &  \leq\lambda\Vert\sum a_{i}e_{k_{j_{i}}}%
\Vert\\
&  \leq\lambda\Vert\sum a_{i}e_{m_{j_{i}}}\Vert\\
&  \leq\lambda|||\sum a_{i}e_{m_{j_{i}}}|||.
\end{align*}
Thus $|||\cdot|||$ is equivalent to $\Vert\cdot\Vert$ on the subspace
$\operatorname*{span}\{(e_{k_{j}})\}.$
\end{proof}

\section{Essentially Finitely Generated Spaces}

The fact that $T[(\mathcal{S}_{n_{i}},\theta_{i})_{i=1}^{k}]$ is isomorphic to
$T_{M}[(\mathcal{S}_{n_{i}},\theta_{i})_{i=1}^{k}]$ was stated by Manoussakis
in \cite{M2}. We present a nominally more general result here. Let us note
that Lopez-Abad and Manoussakis \cite{LAM} has undertaken a thorough study of
mixed Tsirelson spaces generated by finitely many terms.

We shall compute the norm of an element in $T[(\mathcal{S}_{n},\theta
_{n})_{n=1}^{\infty}]$, respectively, $T_{M}[(\mathcal{S}_{n},\theta
_{n})_{n=1}^{\infty}]$, with the help of norming trees. This is derived from
the implicit description of the norms given in equations (\ref{equation 0.1})
and (\ref{equation 0.2}) and have been used in \cite{B, LT2, OT}. An
($(\mathcal{S}_{n})$)-\textbf{admissible tree }(respectively,
\textbf{allowable tree}) is a finite collection of elements $(E_{i}^{m})$,
$0\leq m\leq r,$ $1\leq i\leq k(m)$, in $[{\mathbb{N}}]^{<\infty}$ with the
following properties.

\begin{description}
\item[(i)] $k(0)=1$,

\item[(ii)] For each $m$, $E_{1}^{m}<E_{2}^{m}<\dots<E_{k(m)}^{m}$,

\item[(iii)] Every $E_{i}^{m+1}$ is a subset of some $E_{j}^{m}$,

\item[(iv)] For each $j$ and $m$, the collection $\{E_{i}^{m+1}:E_{i}%
^{m+1}\subseteq E_{j}^{m}\}$ is ${\mathcal{S}}_{n}$-admissible (${\mathcal{S}%
}_{n}$-allowable) for some $n$.
\end{description}

The set $E_{1}^{0}$ is called the \textbf{root} of the tree. The elements
$E_{i}^{m}$ are called \textbf{nodes} of the tree. Given a node $E_{i}^{m},$
$h\left(  E_{i}^{m}\right)  =m$ is called the \textbf{height }of the node
$E_{i}^{m}.$ The height of a tree $\mathcal{T}$ is defined by $H\left(
\mathcal{T}\right)  =\max\{ h\left(  E\right)  :E\in\mathcal{T}\}.$ If
$E_{i}^{n}\subseteq E_{j}^{m}$ and $n>m$, we say that $E_{i}^{n}$ is a
\textbf{descendant} of $E_{j}^{m}$ and $E_{j}^{m}$ is an \textbf{ancestor} of
$E_{i}^{n}$. If, in the above notation, $n=m+1$, then $E_{i}^{n}$ is said to
be an \textbf{immediate successor} of $E_{j}^{m}$, and $E_{j}^{m}$ the
\textbf{immediate predecessor} of $E_{i}^{n}$. Nodes with no descendants are
called \textbf{terminal nodes} or \textbf{leaves} of the tree. We denote the
set of all leaves of a tree $\mathcal{T}$ by $\mathcal{L}( \mathcal{T}).$
Nodes that attain maximal height are called \textbf{base nodes.}

Assign \textbf{tags} to the individual nodes inductively as follows. Let
$t(E_{1}^{0})=1$. If $t(E_{i}^{m})$ has been defined and the collection
$(E_{j}^{m+1})$ of all immediate successors of $E_{i}^{m}$ forms an
${\mathcal{S}}_{k}$-admissible (${\mathcal{S}}_{k}$-allowable) collection,
then define $t(E_{j}^{m+1})=\theta_{k}t(E_{i}^{m})$ for all immediate
successors $E_{j}^{m+1}$ of $E_{i}^{m}.$ If $x\in c_{00}$ and ${\mathcal{T}}$
is an admissible (allowable) tree, let ${\mathcal{T}}x=\sum t(E)\Vert
Ex\Vert_{c_{0}}$ where the sum is taken over all leaves in ${\mathcal{T}}$. It
follows from the implicit description of the norm in $T\left[  \left(
\mathcal{S}_{n},\theta_{n}\right)  _{n=1}^{\infty}\right]  $ (respectively,
$T_{M}\left[  \left(  \mathcal{S}_{n},\theta_{n}\right)  _{n=1}^{\infty
}\right]  $) that $\Vert x\Vert=\max{\mathcal{T}}x$ (respectively,
$|||x|||=\max{\mathcal{T}}x$), with the maximum taken over the set of all
admissible (respectively, allowable) trees. Given a node $E\in\mathcal{T}$
with tag $t\left(  E\right)  =\prod_{i=1}^{m}\theta_{n_{i}},$ define
$o_{\mathcal{T}}\left(  E\right)  =\sum_{i=1}^{m}n_{i}$. When there is no
confusion, we write $o( E) $ instead of $o_{\mathcal{T}}(E) .$

To simplify notation, we shall henceforth denote the spaces $T\left[  \left(
\mathcal{S}_{n},\theta_{n}\right)  _{n=1}^{\infty}\right]  $ and $T_{M}\left[
\left(  \mathcal{S}_{n},\theta_{n}\right)  _{n=1}^{\infty}\right]  $ by $X$
and $X_{M}$ respectively. The norms on these spaces will be denoted by
$\|\cdot\|$ and $\|\cdot\|_{X_{M}}$ respectively.

For a fixed $N\in{\mathbb{N}}$, an ${\mathcal{S}}_{N}$-admissible (-allowable)
tree is a tree satisfying conditions \textbf{(i)}--\textbf{(iii)} above and

\begin{description}
\item[(iv$^{\prime}$)] For each $j$ and $m$, the collection $\{E_{i}%
^{m+1}:E_{i}^{m+1}\subseteq E_{j}^{m}\}$ is ${\mathcal{S}}_{N}$-admissible (-allowable).
\end{description}

It is well known that an ${\mathcal{S}}_{m}$-admissible collection of
${\mathcal{S}}_{n}$-admissible sets is ${\mathcal{S}}_{m+n}$-admissible. The
corresponding fact for the \textquotedblleft allowable\textquotedblright\ case
comes from \cite{ADKM} (see also \cite[Lemma 2.1]{M}).

\begin{lemma}
\label{Lemma1}Given an $\left(  \mathcal{S}_{n}\right)  _{n=1}^{\infty}%
$-admissible (-allowable) tree $\mathcal{T}$ of finite height, there exists an
$\mathcal{S}_{1}$-admissible (-allowable) tree $\mathcal{T}^{\prime}$ with the
same root such that $\mathcal{L}\left(  \mathcal{T}\right)  =\mathcal{L}%
\left(  \mathcal{T}^{\prime}\right)  ,$ and $o_{\mathcal{T}}\left(  E\right)
=o_{\mathcal{T}^{\prime}}\left(  E\right)  $ for all $E\in\mathcal{L}\left(
\mathcal{T}\right)  .$
\end{lemma}

\begin{proof}
The proof is by induction on the height $H\left(  \mathcal{T}\right)  $ of
$\mathcal{T}$. If $H\left(  \mathcal{T}\right)  =0$ then there is nothing to
prove. Assume the statement holds if $H(\mathcal{T}) \leq N$ for some $N.$ Let
$\mathcal{T}$ be an $\left(  \mathcal{S}_{n}\right)  _{n=1}^{\infty}%
$-admissible (-allowable) tree with $H\left(  \mathcal{T}\right)  =N+1.$ Let
$\mathcal{E}_{1}$ be the collection of all nodes of $\mathcal{T}$ at height
$1.$ There exists $n_{0}$ such that $\mathcal{E}_{1}$ is $\mathcal{S}_{n_{0}}%
$-admissible (-allowable). It is easy to see that there is an ${\mathcal{S}%
}_{1}$-admissible (-allowable) tree $\mathcal{T}_{1}$ having the same root as
$\mathcal{T}$ and of height $n_{0}$ such that $\mathcal{L}\left(
\mathcal{T}_{1}\right)  =\mathcal{E}_{1}$ and that every $E\in\mathcal{E}_{1}$
is a leaf of $\mathcal{T}_{1}$ at height $n_{0}.$ If $E\in\mathcal{E}$, then
$\mathcal{T}_{E}=\left\{  F\in\mathcal{T}:F\subseteq E\right\}  $ is an
$\left(  \mathcal{S}_{n}\right)  _{n=1}^{\infty}$-admissible (-allowable) tree
with $H\left(  \mathcal{T}_{E}\right)  \leq N.$ By the inductive hypothesis,
for each $E\in\mathcal{E}_{1},$ there exists an $\mathcal{S}_{1}$-admissible
(-allowable) tree $\mathcal{T}_{E}^{\prime}$ with root $E$ such that
$\mathcal{L}\left(  \mathcal{T}_{E}\right)  =\mathcal{L}\left(  \mathcal{T}%
_{E}^{\prime}\right)  $ and $o_{\mathcal{T}_{E}}\left(  F\right)
=o_{\mathcal{T}_{E}^{\prime}}\left(  F\right)  $ for all $F\in\mathcal{L}%
\left(  \mathcal{T}_{E}\right)  .$

Consider $\mathcal{T}^{\prime}=\mathcal{T}_{1}\cup%
{\displaystyle\bigcup\limits_{E\in\mathcal{E}_{1}}}
\mathcal{T}_{E}^{\prime}.$ Then $\mathcal{T}^{\prime}$ is an $\mathcal{S}_{1}%
$-admissible (-allowable) tree with the same root as $\mathcal{T}$. If
$F\in\mathcal{L}\left(  \mathcal{T}\right)  ,$ then $F\subseteq E$ for some
$E\in\mathcal{E}_{1}$ (since the root cannot be a leaf in this case because
$H\left(  \mathcal{T}\right)  \geq N+1\geq1$). Now $o_{\mathcal{T}}\left(
F\right)  =o_{\mathcal{T}_{E}}\left(  F\right)  + n_{0}$ and $F\in
\mathcal{L}\left(  \mathcal{T}_{E}\right)  .$ Hence $F\in\mathcal{L}\left(
\mathcal{T}_{E}^{\prime}\right)  $ and $o_{\mathcal{T}^{\prime}}\left(
F\right)  =o_{\mathcal{T}_{E}^{\prime}}\left(  F\right)  +n_{0}=o_{\mathcal{T}%
_{E}}\left(  F\right)  +n_{0}=o_{\mathcal{T}}\left(  F\right)  .$ Conversely,
if $F\in\mathcal{L}\left(  \mathcal{T}^{\prime}\right)  ,$ then $F\in
\mathcal{L}\left(  \mathcal{T}_{E}^{\prime}\right)  $ for some $E\in
\mathcal{E}_{1}.$ Thus $F\in\mathcal{L}\left(  \mathcal{T}_{E}\right)  $ and
hence $F\in\mathcal{L}\left(  \mathcal{T}\right)  .$
\end{proof}

\begin{lemma}
\label{Tree}\label{Lemma2} Let $\mathcal{T}$ be an $(\mathcal{S}_{n}%
)_{n=1}^{\infty}$-admissible (-allowable) tree. If $\mathcal{E}$ is a
collection of pairwise disjoint nodes of $\mathcal{T}$ such that $o(E) \leq m$
for all $E\in\mathcal{E}$, then $\mathcal{E}$ is $\mathcal{S}_{m}$-admissible (allowable).
\end{lemma}

\begin{proof}
The proof is by induction on $m.$ The case $m=0$ is clear. Now suppose the
Lemma holds for all $k<m$, $m \geq1$. If the root of $\mathcal{T}$ belongs to
$\mathcal{E}$, then it is the only node in $\mathcal{E}$ and the Lemma clearly
holds. Otherwise, let $k$ be such that the nodes $G_{1} < \dots<G_{q}$ in
$\mathcal{T}$ with height $1$ is $\mathcal{S}_{k}$-admissible (-allowable).
Since each $E \in\mathcal{E}$ is either equal to or is a descendant of some
$G_{i}$, $m \geq o(E) \geq o(G_{i}) = k$. If $m = k$, then $\mathcal{E}
\subseteq\{G_{1},\dots, G_{q}\}$ and thus is ${\mathcal{S}}_{m}$-admissible
(-allowable). If $k<m$, then for each $i$, the subtree $\mathcal{T}_{i}$ with
root $G_{i}$ is an admissible (-allowable) tree such that $o_{\mathcal{T}_{i}%
}(E) \leq m-k$ for all $E \in\mathcal{E} \cap\mathcal{T}_{i}$. By induction,
$E \in\mathcal{E} \cap\mathcal{T}_{i}$ is $\mathcal{S}_{m-k}$-admissible
(-allowable). Therefore, $\mathcal{E}$ is an $\mathcal{S}_{k}$-admissible
(-allowable) collection of $\mathcal{S}_{m-k}$-admissible (-allowable) sets,
and hence an $\mathcal{S}_{m}$-admissible (-allowable) set.
\end{proof}

Given $k\in\mathbb{N}$, let $\left\lceil k\right\rceil $ denote the least
integer greater than or equal to $k.$

\begin{lemma}
\label{Lemma3}Let $\mathcal{T}$ be an $\mathcal{S}_{1}$-admissible
(-allowable) tree. For any $N\in\mathbb{N}$ there exists an $\mathcal{S}_{N}%
$-admissible (-allowable) tree $\mathcal{T}^{\prime}$ with the same root such
that $\mathcal{L}\left(  \mathcal{T}\right)  =\mathcal{L}\left(
\mathcal{T}^{\prime}\right)  $ and $o_{\mathcal{T}^{\prime}}\left(  E\right)
=N\left\lceil o_{\mathcal{T}}\left(  E\right)  /N\right\rceil $ for all
$E\in\mathcal{L}\left(  \mathcal{T}\right)  .$
\end{lemma}

\begin{proof}
Note that the statement holds if $H\left(  \mathcal{T}\right)  \leq N$ by
Lemma \ref{Lemma2}. Now suppose that the statement holds if $H\left(
\mathcal{T}\right)  \leq kN$ for some $k\in\mathbb{N}$. Let $\mathcal{T}$ be
an ${\mathcal{S}}_{1}$-admissible (-allowable) tree with $H\left(
\mathcal{T}\right)  \leq\left(  k+1\right)  N.$ Denote by $\mathcal{T}_{0}$
the tree consisting of all nodes in $\mathcal{T}$ with height $\leq N.$ For
each $E\in\mathcal{T}$ at height $N,$ $H\left(  \mathcal{T}_{E}\right)  \leq
kN$, where $\mathcal{T}_{E}$ consists of all nodes $F$ in $\mathcal{T}$ such
that $F \subseteq E$. By induction, for each $E\in\mathcal{T}$ at height $N$,
there exists an $\mathcal{S}_{N}$-admissible (-allowable) tree $\mathcal{T}%
_{E}^{\prime}$ with root $E$ such that $\mathcal{L}\left(  \mathcal{T}%
_{E}\right)  =\mathcal{L}\left(  \mathcal{T}_{E}^{\prime}\right)  $ and
$o_{\mathcal{T}_{E}^{\prime}}\left(  F\right)  =N\left\lceil o_{\mathcal{T}%
_{E}}\left(  F\right)  /N\right\rceil $ for all $F\in\mathcal{L}\left(
\mathcal{T}_{E}\right)  .$ At the same time, there exists an $\mathcal{S}_{N}%
$-admissible (-allowable) tree $\mathcal{T}_{0}^{\prime}$ with the same root
as $\mathcal{T}_{0}$ such that $\mathcal{L}\left(  \mathcal{T}_{0}^{\prime
}\right)  =\mathcal{L}\left(  \mathcal{T}_{0}\right)  $ and $o_{\mathcal{T}%
^{\prime}_{0}}\left(  F\right)  =N\left\lceil o_{\mathcal{T}_{0}}\left(
F\right)  /N\right\rceil $ for all $F\in\mathcal{L}\left(  \mathcal{T}%
_{0}\right)  .$ Let $\mathcal{T}^{\prime}=\mathcal{T}_{0}^{\prime}\cup%
{\displaystyle\bigcup}
\mathcal{T}_{E}^{\prime},$ where the second union is taken over all nodes
$E\in\mathcal{T}$ at height $N.$ Then $\mathcal{T}^{\prime}$ is an
$\mathcal{S}_{N}$-admissible (-allowable) tree with the same root as
$\mathcal{T}$.

If $E\in\mathcal{L}\left(  \mathcal{T}\right)  $ and $h(E) < N$, then
$E\in\mathcal{L}(\mathcal{T}_{0}) = \mathcal{L}(\mathcal{T}^{\prime}_{0})$ and
has no descendants in $\mathcal{T}^{\prime}$. Hence $E \in\mathcal{L}%
(\mathcal{T}^{\prime})$. Moreover, $o_{\mathcal{T}^{\prime}}\left(  E\right)
=o_{\mathcal{T}_{0}^{\prime}}\left(  E\right)  =N\left\lceil o_{\mathcal{T}%
_{0}}\left(  E\right)  /N\right\rceil =N\left\lceil o_{\mathcal{T}}\left(
E\right)  /N\right\rceil .$ If $E\in\mathcal{L}\left(  \mathcal{T}\right)  $
and $h(E) \geq N,$ then $E\subseteq F$ for some $F\in\mathcal{T}$ at height
$N.$ Hence $E\in\mathcal{L}\left(  \mathcal{T}_{F}\right)  = \mathcal{L}%
\left(  \mathcal{T}_{F}^{\prime}\right)  \subseteq\mathcal{L}\left(
\mathcal{T}^{\prime}\right)  $ and
\begin{align*}
o_{\mathcal{T}^{\prime}}\left(  E\right)   &  = N+o_{\mathcal{T}_{F}^{\prime}%
}\left(  E\right)  =N+N\left\lceil o_{\mathcal{T}_{F}}\left(  E\right)
/N\right\rceil \\
&  =N\left\lceil \frac{o_{\mathcal{T}_{F}}\left(  E\right)  +N}{N}\right\rceil
=N\left\lceil o_{\mathcal{T}}\left(  E\right)  /N\right\rceil .
\end{align*}

Conversely, suppose that $E\in\mathcal{L}(\mathcal{T}^{\prime})$. Then either
$E\in\mathcal{L}\left(  \mathcal{T}_{0}^{\prime}\right)  =\mathcal{L}\left(
\mathcal{T}_{0}\right)  $ with $h(E) < N$ (taken in $\mathcal{T}_{0}$) or else
$E\in\mathcal{L}\left(  \mathcal{T}_{F}^{\prime}\right)  $ for some
$F\in\mathcal{T}$ at height $N.$ Thus $E\in\mathcal{L}\left(  \mathcal{T}%
_{F}\right)  $. In either case, $E \in\mathcal{L}(\mathcal{T})$.
\end{proof}

Combining Lemmas \ref{Lemma1} and \ref{Lemma3}, we obtain:

\begin{proposition}
\label{Prop4}Let $\mathcal{T}$ be an $\left(  \mathcal{S}_{n}\right)
_{n=1}^{\infty}$-admissible (-allowable) tree $\mathcal{T}$ and let
$N\in\mathbb{N}$. Then there exists an $\mathcal{S}_{N}$-admissible
(-allowable) tree $\mathcal{T}^{\prime}$ with the same root such that
$\mathcal{L}\left(  \mathcal{T}\right)  =\mathcal{L}\left(  \mathcal{T}%
^{\prime}\right)  $ and $o_{\mathcal{T}^{\prime}}\left(  E\right)
=N\left\lceil o_{\mathcal{T}}\left(  E\right)  /N\right\rceil $ for all
$E\in\mathcal{L}\left(  \mathcal{T}\right)  .$
\end{proposition}

\begin{proposition}
\label{Prop5}Let $\left(  \theta_{n}\right)  _{n=1}^{\infty}$ be a regular
sequence. Suppose that there exists $N\in\mathbb{N}$ such that $\theta
_{N}^{1/N}=\theta=\sup\theta_{n}^{1/n},$ then the spaces $X, X_{M}, Y$, and
$Y_{M}$ are pairwise isomorphic via the formal identity, where $Y$ and $Y_{M}$
denote the spaces $T[{\mathcal{S}}_{1},\theta]$ and $T_{M}[{\mathcal{S}}%
_{1},\theta]$ respectively.
\end{proposition}

\begin{proof}
It is known that $Y$ and $Y_{M}$ are isomorphic via the formal identity
\cite{B,CO,M}. We shall show that $X_{M}$ to $Y_{M}$ via the formal identity.
The proof that $X$ is isomorphic to $Y$ via the formal identity is similar.
Let $x$ be a finitely supported vector. There exists an $\left(
\mathcal{S}_{n}\right)  _{n=1}^{\infty}$-allowable tree $\mathcal{T}$ such
that
\[
\left\Vert x\right\Vert _{X_{M}}=\sum_{E\in\mathcal{L}\left(  \mathcal{T}%
\right)  }t\left(  E\right)  \left\Vert Ex\right\Vert _{c_{0}}.
\]
By Proposition \ref{Prop4}, there exists an $\mathcal{S}_{1}$-allowable tree
$\mathcal{T}^{\prime}$ with the same root such that $\mathcal{L}\left(
\mathcal{T}\right)  =\mathcal{L}\left(  \mathcal{T}^{\prime}\right)  $ and
$o_{\mathcal{T}^{\prime}}\left(  E\right)  =N\left\lceil o_{\mathcal{T}%
}\left(  E\right)  /N\right\rceil $ for all $E\in\mathcal{L}\left(
\mathcal{T}\right)  .$ If $E\in\mathcal{L}\left(  \mathcal{T}\right)  $ and
$t\left(  E\right)  =\theta_{n_{1}}\cdots\theta_{n_{j}},$ then
\[
t\left(  E\right)  \leq\theta^{n_{1}}\cdots\theta^{n_{j}}=\theta^{n_{1}%
+\dots+n_{j}}=\theta^{o_{\mathcal{T}}\left(  E\right)  }<\theta^{-N}%
\theta^{o_{\mathcal{T}^{\prime}}\left(  E\right)  }.
\]
Therefore,%
\begin{align*}
\left\Vert x\right\Vert _{X_{M}}  &  =\sum_{E\in\mathcal{L}\left(
\mathcal{T}\right)  }t\left(  E\right)  \left\Vert Ex\right\Vert _{c_{0}}\\
&  <\theta^{-N}\sum_{E\in\mathcal{L}\left(  \mathcal{T}^{\prime}\right)
}\theta^{o_{\mathcal{T}^{\prime}}\left(  E\right)  }\left\Vert Ex\right\Vert
_{c_{0}}\leq\theta^{-N}\left\Vert x\right\Vert _{Y_{M}}.
\end{align*}

Conversely, choose an $\mathcal{S}_{1}$-allowable tree $\mathcal{T}%
^{\prime\prime}$ such that
\[
\left\Vert x\right\Vert _{Y_{M}}=\sum_{E\in\mathcal{L}\left(  \mathcal{T}%
^{\prime\prime}\right)  }t\left(  E\right)  \left\Vert Ex\right\Vert _{c_{0}%
}.
\]
Since $\mathcal{T}^{\prime\prime}$ is also $(S_{n})_{n=1}^{\infty}$-allowable,
there exists an $S_{N}$-allowable tree $\mathcal{T}^{\prime\prime\prime}$ such
that $\mathcal{L}\left(  \mathcal{T}^{\prime\prime}\right)  =\mathcal{L}%
\left(  \mathcal{T}^{\prime\prime\prime}\right)  $ and $o_{\mathcal{T}%
^{\prime\prime\prime}}\left(  E\right)  =N\left\lceil o_{\mathcal{T}%
^{\prime\prime}}\left(  E\right)  /N\right\rceil $ for all $E\in
\mathcal{L}\left(  \mathcal{T}^{\prime\prime}\right)  .$ Hence $t(E)=\theta
^{o_{\mathcal{T}^{\prime\prime}}\left(  E\right)  }\leq\theta
^{-N+o_{\mathcal{T}^{\prime\prime\prime}}\left(  E\right)  }.$ Thus
\[
\left\Vert x\right\Vert _{Y_{M}}\leq\theta^{-N}\sum_{E\in\mathcal{L}\left(
\mathcal{T}^{\prime\prime\prime}\right)  }\theta_{N}^{o_{\mathcal{T}%
^{\prime\prime\prime}}\left(  E\right)  /N}\left\Vert Ex\right\Vert _{c_{0}%
}\leq\left\Vert x\right\Vert _{X_{M}}.
\]
The final inequality holds since $\mathcal{T}^{\prime\prime\prime}$ is also
$(S_{n})_{n=1}^{\infty}$-allowable and the tag of $E$ in $\mathcal{T}%
^{\prime\prime\prime}$ is $\theta_{N}^{o_{\mathcal{T}^{\prime\prime\prime}%
}(E)/N}$.
\end{proof}

\section{Main Construction}

The main aim of the present paper is to show that the spaces $X$ and $X_{M}$
are not isomorphic for a large class or regular sequences $(\theta_{n})$. In
view of Proposition \ref{Nonisomorphic}, it suffices to show that the norms
$\|\cdot\|$ and $\|\cdot\|_{X_{M}}$ are not equivalent on
$\operatorname*{span}\{e_{k}:{k\in V}\}$ for any $V \in[{\mathbb{N}}]$. Our
strategy is to construct, for any $V \in[{\mathbb{N}}]$, vectors $x
\in\operatorname*{span}\{e_{k}:{k\in V}\}$ where the ratio $\|x\|_{X_{M}%
}/\|x\|$ can be made arbitrarily large. The basic units of the construction
are the \textbf{repeated averages} due to Argyros, Mercourakis and Tsarpalias
\cite{A}. These are then layered together, where each layer consists of
repeated averages whose complexities go through a cycle. This variation
\emph{within} a layer is the main feature that distinguishes the present
construction from related previous constructions that are used in, e.g.,
\cite{ADKM, LT2}. The reason for layered construction of vectors is to dictate
that the norming trees that approximately norm the given vector must
structurally resemble the vector itself. In the presence of a condition such
as ($\dag$), one may exploit the large ratio between $\theta_{m+n}$ and
$\theta_{m}\theta_{n}$ to ensure that different layers behave differently. In
the absence of such a condition, one must find a way to ``lock in" the
behavior of the norming tree on the given vector. Our idea is to make the
vector cycle through different complexities within each layer so that the
norming tree is forced to follow these ups and downs.

If $x, y \in\operatorname*{span}\{(e_{k})\}$, we define $x < y$, respectively,
$x \subseteq y$, to mean $\operatorname*{supp} x < \operatorname*{supp} y$ and
$\operatorname*{supp} x \subseteq\operatorname*{supp} y$, respectively. We
shall also say that $E \subseteq x$ if $E \in[{\mathbb{N}}]^{<\infty}$ and $E
\subseteq\operatorname*{supp} x$. An $\mathcal{S}_{0}$-repeated average is a
vector $e_{k}$ for some $k\in\mathbb{N}.$ For any $p\in\mathbb{N}$, an
$\mathcal{S}_{p}$-repeated average is a vector of the form $\frac{1}{k}%
\sum_{i=1}^{k}x_{i},$ where $x_{1}<...<x_{k}$ are repeated $\mathcal{S}_{p-1}%
$-repeated averages and $k=\min\operatorname*{supp}x_{1}.$ Observe that any
$\mathcal{S}_{p}$-repeated average $x$ is a convex combination of $\left\{
e_{k}:k\in\operatorname*{supp}x\right\}  $ such that $\left\Vert x\right\Vert
_{\infty}\leq\left(  \min\operatorname*{supp}x\right)  ^{-1}$ and
$\operatorname*{supp} x \in{\mathcal{S}}_{p}$.

Let $\left(  \theta_{n}\right)  _{n=1}^{\infty}$ be a given regular decreasing
sequence that satisfies the following:

\begin{description}
\item[$\left(  \lnot\dag\right)  $] \label{Nodagger}$\lim_{m}\delta_{m}=0,$
where $\delta_{m}=\limsup_{n}\frac{\theta_{m+n}}{\theta_{n}}.$

\item[$(\ddag)$] \label{doubledagger}There exists $F:\mathbb{N\rightarrow R}$
with $\lim_{n\rightarrow\infty}F\left(  n\right)  =0$ such that for all $R,
t\in\mathbb{N}$ and any arithmetic progression $\left(  s_{i}\right)
_{i=1}^{R}$ in $\mathbb{N}$,
\[
\max_{1\leq i\leq R}\frac{\theta_{s_{i}+t}}{\theta_{s_{i}}}\leq F\left(
R\right)  \sum_{i=1}^{R}\frac{\theta_{s_{i}+t}}{\theta_{s_{i}}}.
\]

\end{description}

\noindent Recall from \S 2 that $X$ and $X_{M}$ are known to be non-isomorphic
if condition ($\dag$) holds. The condition ($\ddag$) is imposed to make the
construction work. As we shall see, it is general enough to include many
interesting cases.

From here on fix $N\in\mathbb{N}$ and $V\in\left[  \mathbb{N}\right]  $
arbitrarily$.$ Choose sequences $\left(  p_{k}\right)  _{k=1}^{N}$ and
$\left(  L_{k}\right)  _{k=1}^{N}$ in $\mathbb{N}$, $L_{k} \geq2$, that
satisfy the following conditions:

\begin{description}
\item[(A)] $\dfrac{\theta_{p_{M+1}+n}}{\theta_{n}}\leq\frac{\theta_{1}%
}{24N^{2}}\prod_{i=1}^{M}\theta_{L_{i}p_{i}}$ if $0\leq M\leq N-2$ and $n\geq
p_{N}$ (the vacuous product $\prod_{i=1}^{0}\theta_{L_{i}p_{i}}$ is taken to
be $1$),

\item[(B)] $p_{M+1}>\sum_{i=1}^{M}L_{i}p_{i}$ if $0 < M \leq N-2$,

\item[(C)] $F\left(  L_{M+1}\right)  \leq\frac{\theta_{1}}{144N^{2}}%
\prod_{i=1}^{M}\theta_{L_{i}p_{i}}$ if $0<M\leq N-2.$
\end{description}

\noindent Note that condition \textbf{(A)} may be realized because of
($\lnot\dag$) and condition \textbf{(C)} by way of (\ddag). Given
$k\in\mathbb{N}$ and $1\leq M\leq N,$ define $r_{M}\left(  k\right)  $ to be
the integer in $\left\{  1,2,...,L_{M}\right\}  $ such that $L_{M}%
|(k-r_{M}\left(  k\right)  ).$ We can construct sequences of vectors
$\mathbf{x}^{0},\mathbf{x}^{1},\dots,\mathbf{x}^{N}$ with the following properties.

\begin{description}
\item[$\boldsymbol{(\alpha)}$] $\mathbf{x}^{0}$ is a subsequence of $\left(
e_{k}\right)  _{k\in V}.$

\item[$\boldsymbol{(\beta)}$] Say $\mathbf{x}^{M}=(x_{j}^{M})$ and $m_{j}%
=\min\operatorname*{supp}x_{j}^{M}.$ Then there is a sequence $(I_{k}^{M+1})$
of integer intervals such that $I_{k}^{M+1}<I_{k+1}^{M+1}$, $%
{\displaystyle\bigcup\limits_{k=1}^{\infty}}
I_{k}^{M+1}=\mathbb{N}$ and each vector $x_{k}^{M+1}\in\mathbf{x}^{M+1}$ is of
the form%
\[
x_{k}^{M+1}=\sum_{j\in I_{k}^{M+1}}a_{j}x_{j}^{M},
\]
where $\theta_{r_{M+1}(k)p_{M+1}}\sum_{j\in I_{k}^{M+1}}a_{j}e_{m_{j}}$ is an
$\mathcal{S}_{r_{M+1}(k)p_{M+1}}$-repeated average. Moreover, the sequence
$(a_{j})_{j=1}^{\infty}$ is decreasing.
\end{description}

Each $x_{k}^{M+1}$ is made up of components of diverse complexities. In order
to estimate its $\Vert\cdot\Vert$- and $\Vert\cdot\Vert_{X_{M}}$- norms, we
decompose $x_{k}^{M+1}$ into components of pure forms in the following manner.
The coefficients $\left(  a_{j}\right)  $ are as given in $\boldsymbol{(\beta
)}$.\newline

\begin{notation}
{Given }$1\leq r_{i}\leq L_{i}${, }$1\leq M\leq N-1,${ write }%
\[
x_{k}^{M+1}\left(  r_{M}\right)  =\sum\limits_{\substack{j\in I_{k}%
^{M+1}\\r_{M}\left(  j\right)  =r_{M}}}a_{j}x_{j}^{M}.
\]
{For }$1\leq s<M,${ define}$\ \ ${\ }%
\[
x_{k}^{M+1}\left(  r_{s},...,r_{M}\right)  =%
{\displaystyle\sum\limits_{\substack{j\in I_{k}^{M+1}\\r_{M}\left(  j\right)
\ =r_{M}}}}
a_{j}x_{j}^{M}\left(  r_{s},...,r_{M-1}\right)  .
\]
{If }$1\leq s\leq M,$ {it is clear that }$x_{k}^{M+1}=\sum x_{k}^{M+1}\left(
r_{s},...,r_{M}\right)  ,${ where the sum is taken over all possible values of
}$r_{s},...,r_{M}.$
\end{notation}

Given a sequence $\mathbf{u}=\left(  u_{1},u_{2},...\right)  $ of linearly
independent vectors, write $\left[  y\right]  _{\mathbf{u}}=\left(
a_{k}\right)  $ if $y=\sum a_{k}u_{k}$. For instance, $\|[x_{k}^{M+1}%
]_{_{\mathbf{x}^{M}}}\|_{\ell^{1}}=\sum_{j\in I_{k}^{M+1}}a_{j}=\theta
_{r_{M+1}\left(  k\right)  p_{M+1}}^{-1}.$ To compute $\|[ x_{k}^{M+1}]
_{_{\mathbf{x}^{s-1}}}]\|_{\ell^{1}}, 1\leq s\leq M,$ calculate the $\ell^{1}%
$-norms of each of the pure forms $[ x_{k}^{M+1}( r_{s},...,r_{M})]
_{\mathbf{x}^{s-1}}$ and sum over all $r_{s},...,r_{M}$.

The following simple Lemma is useful for our computations. A subset $I$ of
${\mathbb{N}}$ is said to be $L$-\textbf{skipped} if $|i-j| \geq L$ whenever
$i$ and $j$ are distinct elements of $I$.

\begin{lemma}
\label{Lskippedseq}If $\left(  a_{i}\right)  $ is a nonnegative decreasing
sequence defined on an interval $J$ in $\mathbb{N}$ and $I$ is an $L$-skipped
set, then
\[
\sum_{i\in I}a_{i}\leq\frac{1}{L}\sum a_{i}+\sup a_{i}.
\]
Moreover, if there exists $r$ such that $I = \{i \in J: i = r \mod L\}$, then
\[
\frac{1}{L}\sum a_{i}-\sup a_{i}\leq\sum_{i\in I}a_{i}.
\]

\end{lemma}

\begin{proposition}
\label{ellonenormsrel}If $1 \leq s \leq M < N$ and $k\in{\mathbb{N}},$ then
\[
\prod_{i=s}^{M}\left(  L_{i}^{-1}-k^{-1}\right)  \leq\frac{\left\Vert \left[
x_{k}^{M+1}\left(  r_{s},\cdots,r_{M}\right)  \right]  _{\mathbf{x}^{s-1}%
}\right\Vert _{\ell^{1}}}{\theta_{r_{M+1}\left(  k\right)  p_{M+1}}^{-1}%
\prod_{i=s}^{M}\theta_{r_{i}p_{i}}^{-1}}\leq\prod_{i=s}^{M}\left(  L_{i}%
^{-1}+k^{-1}\right)  .
\]

\end{proposition}

\begin{proof}
The proof is by induction on $M.$ When $M=s,$%
\begin{align*}
\left\Vert \left[  x_{k}^{M+1}\left(  r_{s},...,r_{M}\right)  \right]
_{\mathbf{x}^{s-1}}\right\Vert _{\ell^{1}}  &  =\left\Vert \left[  x_{k}%
^{M+1}\left(  r_{M}\right)  \right]  _{\mathbf{x}^{M-1}}\right\Vert _{\ell
^{1}}\\
&  =\biggl\Vert \Biggl[ \sum\limits_{\substack{j\in I_{k}^{M+1}\\r_{M}(j)
=r_{M}}}a_{j}x_{j}^{M}\Biggr]_{_{\mathbf{x}^{M-1}}}\biggr\Vert _{\ell^{1}}\\
&  =\sum\limits_{\substack{j\in I_{k}^{M+1}\\r_{M}\left(  j\right)  =r_{M}%
}}a_{j}\left\Vert \left[  x_{j}^{M}\right]  _{_{\mathbf{x}^{M-1}}}\right\Vert
_{\ell^{1}}\\
&  =\theta_{r_{M}p_{M}}^{-1}\sum\limits_{\substack{j\in I_{k}^{M+1}%
\\r_{M}\left(  j\right)  =r_{M}}}a_{j}.
\end{align*}
Note that $\left\{  j\in I_{k}^{M+1}:r_{M}\left(  j\right)  =r_{M}\right\}  $
is an $L_{M}$-skipped subset of the integer interval $I_{k}^{M+1}$. It follows
from Lemma \ref{Lskippedseq} that%
\begin{align*}
\sum\limits_{\substack{j\in I_{k}^{M+1}\\r_{M}\left(  j\right)  =r_{M}}}a_{j}
&  \leq\frac{1}{L_{M}}\sum\limits_{j\in I_{k}^{M+1}}a_{j}+\sup a_{j}\\
&  \leq\left(  L_{M}^{-1}+k^{-1}\right)  \theta_{r_{M+1}\left(  k\right)
p_{M+1}}^{-1}.
\end{align*}
Therefore, $\left\Vert \left[  x_{k}^{M+1}\left(  r_{M}\right)  \right]
_{\mathbf{x}^{M-1}}\right\Vert _{\ell^{1}}\leq\theta_{r_{M+1}\left(  k\right)
p_{M+1}}^{-1}\theta_{r_{M}p_{M}}^{-1}\left(  L_{M}^{-1}+k^{-1}\right)  .$

Suppose that the Proposition holds for $M-1$. Then%
\begin{align*}
\Vert[ x_{k}^{M+1}  &  (r_{s},...,r_{M})]_{\mathbf{x}^{s-1}}\Vert_{\ell^{1}} &
\\
&  =\Vert\sum\limits_{\substack{j\in I_{k}^{M+1}\\r_{M}\left(  j\right)
=r_{M}}}a_{j}\left[  x_{j}^{M}\left(  r_{s},...,r_{M-1}\right)  \right]
_{\mathbf{x}^{s-1}}\Vert_{\ell^{1}} & \\
&  =\sum\limits_{\substack{j\in I_{k}^{M+1}\\r_{M}\left(  j\right)  =r_{M}%
}}a_{j}\left\Vert \left[  x_{j}^{M}\left(  r_{s},...,r_{M-1}\right)  \right]
_{\mathbf{x}^{s-1}}\right\Vert _{\ell^{1}} & \\
&  \leq\prod_{i=s}^{M-1}\theta_{r_{i}p_{i}}^{-1}\left(  L_{i}^{-1}%
+k^{-1}\right)  \cdot\sum\limits_{\substack{j\in I_{k}^{M+1}\\r_{M}\left(
j\right)  =r_{M}}}\frac{a_{j}}{\theta_{r_{M}\left(  j\right)  p_{M}}} & \\
&  \quad\quad\text{ by the inductive hypothesis,} & \\
&  =\prod_{i=s}^{M-1}\theta_{r_{i}p_{i}}^{-1}\left(  L_{i}^{-1}+k^{-1}\right)
\cdot\theta_{r_{M}p_{M}}^{-1}\sum\limits_{\substack{j\in I_{k}^{M+1}%
\\r_{M}\left(  j\right)  =r_{M}}}a_{j} & \\
&  \leq\prod_{i=s}^{M-1}\theta_{r_{i}p_{i}}^{-1}\left(  L_{i}^{-1}%
+k^{-1}\right)  \cdot\theta_{r_{M}p_{M}}^{-1}\left(  \frac{1}{L_{M}}%
\sum\limits_{j\in I_{k}^{M+1}}a_{j}+\sup a_{j}\right)  \text{ } & \\
&  \quad\quad\text{ by Lemma \ref{Lskippedseq},} & \\
&  =\prod_{i=s}^{M-1}\theta_{r_{i}p_{i}}^{-1}\left(  L_{i}^{-1}+k^{-1}\right)
\theta_{r_{M}p_{M}}^{-1}\theta_{r_{M+1}\left(  k\right)  p_{M+1}}^{-1}\left(
L_{M}^{-1}+k^{-1}\right)  & \\
&  =\theta_{r_{M+1}\left(  k\right)  p_{M+1}}^{-1}\prod_{i=s}^{M}\theta
_{r_{i}p_{i}}^{-1}\left(  L_{i}^{-1}+k^{-1}\right)  . &
\end{align*}
The other inequality is proved similarly.
\end{proof}

\noindent From this point onwards, we shall only consider those $k$'s that
satisfy
\begin{equation}
k\geq42N^{2}\prod_{i=1}^{N}L_{i}\theta_{L_{i}p_{i}}^{-1}. \label{1}%
\end{equation}
It follows from the choice of $k$ that for all $1\leq s\leq M\leq N,$
\begin{equation}
\prod_{i=s}^{M}\left(  L_{i}^{-1}+k^{-1}\right)  \leq2\prod_{i=s}^{M}%
L_{i}^{-1}. \label{2}%
\end{equation}
Indeed, since $L_{i}^{-1}+k^{-1} \leq\left(  1+\frac{1}{42N}\right)
L^{-1}_{i}$ for all $i,$ we have%
\begin{align*}
\prod_{i=s}^{M}( L_{i}^{-1}+k^{-1})  &  \leq(1+\frac{1}{42N})^{N}\prod
_{i=s}^{M}L_{i}^{-1}\\
&  \leq e^{1/42}\prod_{i=s}^{M}L_{i}^{-1}<2\prod_{i=s}^{M}L_{i}^{-1}.
\end{align*}
Likewise, for all $1\leq s\leq M\leq N,$%
\begin{equation}
\prod_{i=s}^{M}\left(  L_{i}^{-1}-k^{-1}\right)  > \frac{1}{2}\prod_{i=s}%
^{M}L_{i}^{-1}. \label{3}%
\end{equation}

\begin{corollary}
\label{CorNorm1}If $1 \leq s \leq M < N$ and $k$ satisfies (\ref{1}), then
\[
\frac{1}{2}\leq\frac{\left\Vert \left[  x_{k}^{M+1}\left(  r_{s},\cdots
,r_{M}\right)  \right]  _{\mathbf{x}^{s-1}}\right\Vert _{\ell^{1}}}%
{\theta_{r_{M+1}\left(  k\right)  p_{M+1}}^{-1}\prod_{i=s}^{M}\theta
_{r_{i}p_{i}}^{-1}L_{i}^{-1}}\leq2.
\]

\end{corollary}

\begin{corollary}
\label{CorNorm}If $k$ satisfies (\ref{1}) and $1 \leq M \leq N$, then
\[
\left\Vert x_{k}^{M}\right\Vert _{\ell^{1}}\leq2\prod_{i=1}^{M}\theta
_{L_{i}p_{i}}^{-1}.
\]

\end{corollary}

\begin{proof}
If $M=1,$ then
\[
\|x^{1}_{k}\|_{\ell^{1}} = \|[x^{1}_{k}]_{\mathbf{x}^{0}}\|_{\ell^{1}} =
\theta^{-1}_{r_{1}(k)p_{1}} \leq\theta^{-1}_{L_{1}p_{1}}.
\]
If $M\geq2,$ according to Corollary \ref{CorNorm1},
\begin{align*}
\left\Vert x_{k}^{M}\right\Vert _{\ell^{1}}  &  =\sum_{r_{1},\cdots,r_{M-1}%
}\left\Vert \left[  x_{k}^{M}\left(  r_{1},...,r_{M-1}\right)  \right]
_{\mathbf{x}^{0}}\right\Vert _{\ell^{1}}\\
&  \leq2\theta_{r_{M}\left(  k\right)  p_{M}}^{-1}\sum_{r_{1},...,r_{M-1}%
}\prod_{i=1}^{M-1}\theta_{r_{i}p_{i}}^{-1}L_{i}^{-1}\\
&  \leq2\theta_{L_{M}p_{M}}^{-1}\prod_{i=1}^{M-1}\theta_{L_{i}p_{i}}^{-1}
=2\prod_{i=1}^{M}\theta_{L_{i}p_{i}}^{-1}\text{.}%
\end{align*}

\end{proof}

We shall employ the same decomposition technique to estimate $\left\Vert
x_{k}^{N}\right\Vert _{X_{M}}.$ To simplify notation, let $p(r_{M}%
,\dots,r_{M^{\prime}}) = \sum^{M^{\prime}}_{i=M}p_{i}r_{i}$ if $M \leq
M^{\prime}$.

\begin{proposition}
\label{ModNorm} If $k$ satisfies (\ref{1}), then%
\[
\left\Vert x_{k}^{N}\right\Vert _{X_{M}}\geq\frac{\theta_{1}}{2}\sum
_{r_{1},...,r_{N-1}}\theta_{p\left(  r_{1},...,r_{N-1},r_{N}\left(  k\right)
\right)  }\theta_{r_{N}\left(  k\right)  p_{N}}^{-1}\prod_{i=1}^{N-1}%
\theta_{r_{i}p_{i}}^{-1}L_{i}^{-1}.
\]

\end{proposition}

\begin{proof}
We first decompose $x_{k}^{N}$ into a sum of pure forms, i.e.,
\[
x_{k}^{N}=\sum_{r_{1},...,r_{N-1}}x_{k}^{N}\left(  r_{1},...,r_{N-1}\right)
.
\]
Now given $r_{1},...,r_{N-1},$ $\operatorname*{supp}x_{k}^{N}\left(
r_{1},...,r_{N-1}\right)  \in\mathcal{S}_{p\left(  r_{1},...,r_{N-1}%
,r_{N}\left(  k\right)  \right)  }.$ Hence
\begin{align*}
\left\Vert x_{k}^{N}\left(  r_{1},...,r_{N-1}\right)  \right\Vert _{X_{M}}  &
\geq\theta_{p\left(  r_{1},...,r_{N-1},r_{N}\left(  k\right)  \right)
}\left\Vert x_{k}^{N}\left(  r_{1},...,r_{N-1}\right)  \right\Vert _{\ell^{1}%
}\\
&  \geq\frac{\theta_{p(r_{1},...,r_{N-1},r_{N}(k))}}{2}\theta_{r_{N}\left(
k\right)  p_{N}}^{-1}\prod_{i=1}^{N-1}\theta_{r_{i}p_{i}}^{-1}L_{i}^{-1}\text{
}%
\end{align*}
by Corollary \ref{CorNorm}. Since $k\geq\prod_{i=1}^{N-1}L_{i}$ by (\ref{1}),
$S\in$ $\mathcal{S}_{1}$ whenever $S\subseteq\mathbb{N}$ satisfies $k\leq\min
S$ and $\left\vert S\right\vert \leq\prod_{i=1}^{N-1}L_{i}.$ In particular,
\[
\left\{  \operatorname*{supp}x_{k}^{N}\left(  r_{1},...,r_{N-1}\right)  :1\leq
r_{i}\leq L_{i},\text{ }1\leq i\leq N-1\right\}
\]
is $\mathcal{S}_{1}$-allowable. Thus
\begin{align*}
\left\Vert x_{k}^{N}\right\Vert _{X_{M}}  &  \geq\theta_{1}\sum_{r_{1}%
,...,r_{N-1}}\left\Vert x_{k}^{N}\left(  r_{1},...,r_{N-1}\right)  \right\Vert
_{X_{M}}\\
&  \geq\frac{\theta_{1}}{2}\sum_{r_{1},...,r_{N-1}}\theta_{p\left(
r_{1},...,r_{N-1},r_{N}\left(  k\right)  \right)  }\theta_{r_{N}\left(
k\right)  p_{N}}^{-1}\prod_{i=1}^{N-1}\theta_{r_{i}p_{i}}^{-1}L_{i}^{-1},
\end{align*}
as required.
\end{proof}

The following estimate is easily obtainable from Proposition \ref{ModNorm}.

\begin{corollary}
\label{ModNormM}If $0\leq M < N-1$ and $k$ satisfies (\ref{1}), then%
\begin{equation}
\left\Vert x_{k}^{N}\right\Vert _{X_{M}}\geq\frac{\theta_{1}}{2}\sum
_{r_{M+1},...,r_{N-1}}\theta_{p\left(  r_{M+1},...,r_{N-1},r_{N}\left(
k\right)  \right)  }\theta_{r_{N}\left(  k\right)  p_{N}}^{-1}\prod
_{i=M+1}^{N-1}\theta_{r_{i}p_{i}}^{-1}L_{i}^{-1}. \label{B'}%
\end{equation}

\end{corollary}

\begin{proof}
By Proposition \ref{ModNorm} and the regularity of $\left(  \theta_{n}\right)
$,
\begin{align*}
\left\Vert x_{k}^{N}\right\Vert _{X_{M}}  &  \geq\frac{\theta_{1}}{2}%
\sum_{r_{1},...,r_{N-1}}\theta_{p\left(  r_{1},...,r_{N}\left(  k\right)
\right)  }\theta_{r_{N}\left(  k\right)  p_{N}}^{-1}\prod_{i=1}^{N-1}%
\theta_{r_{i}p_{i}}^{-1}L_{i}^{-1}\\
&  \geq\frac{\theta_{1}}{2}\sum_{r_{1},...,r_{N-1}}\theta_{p\left(
r_{2},...,r_{N}\left(  k\right)  \right)  }\theta_{r_{1}p_{1}}\theta
_{r_{N}\left(  k\right)  p_{N}}^{-1}\prod_{i=1}^{N-1}\theta_{r_{i}p_{i}}%
^{-1}L_{i}^{-1}\\
&  =\frac{\theta_{1}}{2}\sum_{r_{2},...,r_{N-1}}\theta_{p\left(
r_{2},...,r_{N}\left(  k\right)  \right)  }\theta_{r_{N}\left(  k\right)
p_{N}}^{-1}\prod_{i=2}^{N-1}\theta_{r_{i}p_{i}}^{-1}L_{i}^{-1}.
\end{align*}
Repeat the argument $M$ times to obtain the required result.
\end{proof}

The main bulk of the calculations occur in estimating the $X$-norm of
$x^{N}_{k}$. The next lemma is the mechanism behind one of the crucial
estimates (Proposition \ref{PropA}). If $x \in c_{00}$ and $p \geq0$, let
$\|x\|_{{\mathcal{S}}_{p}} = \sup_{E\in{\mathcal{S}}_{p}}\|Ex\|_{\ell^{1}}$.

\begin{lemma}
\label{L1pq}Let $p,q\geq0$, and $P=(m_{n})\in\lbrack{\mathbb{N}}]$ be given.
Assume that $G_{1}<G_{2}<\cdots$ is a sequence in $\left[  P\right]
^{<\infty}$ such that $\ \sum_{m_{n}\in G_{i}}a_{n}e_{m_{n}}$ is an
$\mathcal{S}_{q}$-repeated average for all $i$ and that there exists
$Q=(m_{n_{k}})\in\lbrack P]$ so that for each $k$, there is a vector $z_{k}$ satisfying:

\begin{enumerate}
\item $\operatorname{supp}z_{k}\subseteq\lbrack m_{n_{k}},m_{n_{k}+1})$,

\item $\Vert z_{k}\Vert_{\ell^{1}}\leq1$,

\item $\Vert\sum_{k=1}^{j}z_{k}\Vert_{{\mathcal{S}}_{p}}\leq6$ for all
$j\in{\mathbb{N}}$.
\end{enumerate}

\noindent Set $y_{i}=\sum_{m_{n_{k}}\in G_{i}}a_{n_{k}}z_{k}$. Then

\begin{enumerate}
\item[(i)] $\Vert\sum_{i=1}^{j}y_{i}\Vert_{{\mathcal{S}}_{p+q}}\leq6$ for all
$j\in{\mathbb{N}}$,

\item[(ii)] $\Vert y_{i}\Vert_{{\mathcal{S}}_{p+q-1}}\leq6/m$ if $q\geq1$,
where $m=\min G_{i}$.
\end{enumerate}
\end{lemma}

\begin{proof}
We first establish (i). The proof is by induction on $q.$ The case $q=0$ is
trivial. Assume the result holds for some $q,$ we shall prove it for $q+1.$ If
$G_{1}<G_{2}<\cdots$ is a sequence in $\left[  P\right]  ^{<\infty}$ such that
$\sum_{m_{n}\in G_{i}}a_{n}e_{m_{n}}$ is an $\mathcal{S}_{q+1}$-repeated
average for all $i$, then each of these $\mathcal{S}_{q+1}$-repeated averages
can be written as $\frac{1}{m_{n\left(  i\right)  }}\sum_{t\in H_{i}}%
\sum_{m_{n}\in F_{t}}b_{n}e_{m_{n}},$ where $m_{n\left(  i\right)  }=\min
G_{i}=|H_{i}|$, $F_{t}<F_{t^{\prime}}$ if $t<t^{\prime}$ and $\sum_{m_{n}\in
F_{t}}b_{n}e_{m_{n}}$ is an $\mathcal{S}_{q}$-repeated average for all $t.$
Let $y_{i}=\sum_{m_{n_{k}}\in G_{i}}a_{n_{k}}z_{k}.$ Then $y_{i}=\frac
{1}{m_{n\left(  i\right)  }}\sum_{t\in H_{i}}v_{t},$ where $v_{t}%
=\sum_{m_{n_{k}}\in F_{t}}b_{n_{k}}z_{k}.$ Given a set $J\in\mathcal{S}%
_{p+q+1},$ write $J=\cup_{l=1}^{s}J_{l},$ $J_{1}<\cdots<J_{s},$ $J_{l}%
\in\mathcal{S}_{p+q},$ $s\leq\min J.$ Note that by induction, $\left\Vert
J_{l}\left(  \sum_{t\in H_{i}}v_{t}\right)  \right\Vert _{\ell^{1}}\leq6$ for
all $l$ and $i.$ Hence $\left\Vert J_{l}y_{i}\right\Vert _{\ell^{1}}\leq
\frac{6}{m_{n\left(  i\right)  }}.$ Let $i_{0}$ be the smallest number such
that $J\cap\operatorname*{supp}z_{k}\neq\emptyset$ for some $m_{n_{k}}\in
H_{i_{0}}.$ For any $j$,%
\begin{align*}
\left\Vert J(\sum_{i=1}^{j}y_{i})\right\Vert _{\ell^{1}}  &  \leq\sum
_{i=i_{0}}^{i_{0}+2}\left\Vert y_{i}\right\Vert _{\ell^{1}}+\sum_{l}%
\sum_{i=i_{0}+3}^{\infty}\left\Vert J_{l}y_{i}\right\Vert _{\ell^{1}}\\
&  \leq3+\sum_{l}\sum_{i=i_{0}+3}^{\infty}\frac{6}{m_{n\left(  i\right)  }}\\
&  \leq3+\sum_{l}\frac{12}{m_{n\left(  i_{0}+3\right)  }}\text{, \ \ since
}m_{n\left(  i+1\right)  }\geq2m_{n\left(  i\right)  },\\
&  =3+\frac{12s}{m_{n\left(  i_{0}+3\right)  }}.
\end{align*}
But since $\min J<m_{n\left(  i_{0}+1\right)  },$ $s/m_{n\left(
i_{0}+3\right)  }<1/4.$ Therefore,%
\[
\left\Vert J(\sum_{i=1}^{j}y_{i})\right\Vert _{\ell^{1}}<3+\frac{12}{4}=6.
\]
To prove (ii), note that an $\mathcal{S}_{q}$-repeated average $\sum_{m_{n}\in
G_{i}}a_{n}e_{m_{n}}$ may be written as $m^{-1}\left(  u_{1}+\cdots
+u_{m}\right)  ,$ where $u_{1}<\cdots<u_{m}$ are $\mathcal{S}_{q-1}$-repeated
averages. If $u_{j}=\sum_{m_{n}\in F_{j}}b_{n}e_{m_{n}},$ then $y_{i}=$
$m^{-1}\left(  w_{1}+\cdots+w_{m}\right)  ,$ where $w_{j}=\sum_{m_{n_{k}}\in
F_{j}}b_{n_{k}}z_{k}.$ By (i), if $J\in\mathcal{S}_{p+q-1},$ then $\left\Vert
J\left(  w_{1}+\cdots+w_{m}\right)  \right\Vert _{\ell^{1}}\leq6$. Hence
$\left\Vert Jy_{i}\right\Vert _{\ell^{1}}\leq6/m.$
\end{proof}

Assume that $0\leq M<M+s\leq N$ and that $r_{1},..., r_{N-1}$ are given. For
notational convenience, let $x_{k}^{M+s}\left(  r_{M+1},\dots, r_{M+s-1}%
\right)  =x_{k}^{M+s}$ if $s=1.$ Taking $m_{j}=\min\operatorname*{supp}%
x_{j}^{M},$ define
\[
u_{k}^{M+s}\left(  r_{M+1},\dots, r_{M+s-1}\right)  =\sum b_{j}e_{m_{j}}%
\]
if $x_{k}^{M+s}\left(  r_{M+1},\dots, r_{M+s-1}\right)  =\sum b_{j}x_{j}^{M}$.
(The vector is also labeled as $u^{M+1}_{k}$ if $s =1$.)

\begin{proposition}
\label{P1}Let $r_{N}=r_{N}\left(  k\right)  .$ Then
\[
\left\Vert u_{k}^{N}\left(  r_{M+1},\dots, r_{N-1}\right)  \right\Vert
_{\mathcal{S}_{p\left(  r_{M+1},\dots, r_{N}\left(  k\right)  \right)  -1}%
}\leq\frac{6}{k}\prod_{i=M+1}^{N}\theta_{r_{i}p_{i}}^{-1}.
\]

\end{proposition}

\begin{proof}
We shall apply Lemma \ref{L1pq} repeatedly to show that
\begin{equation}
\prod_{i=M+1}^{M+s}\theta_{r_{i}p_{i}}\left\Vert u_{t}^{M+s}\left(
r_{M+1},...,r_{M+s-1}\right)  \right\Vert _{\mathcal{S}_{p\left(
r_{M+1},...,r_{M+s}\right)  -1}}\leq\frac{6}{t} \label{p10}%
\end{equation}
if $r_{M+s}\left(  t\right)  =r_{M+s}$ and%
\[
\prod_{i=M+1}^{M+s}\theta_{r_{i}p_{i}}\left\Vert \sum u_{t}^{M+s}\left(
r_{M+1},...,r_{M+s-1}\right)  \right\Vert _{\mathcal{S}_{p\left(
r_{M+1},...,r_{M+s}\right)  }}\leq6
\]
for any sum over a finite set of $t$'s satisfying $r_{M+s}\left(  t\right)
=r_{M+s}.$ Suppose that $s=1.$ Set $p=0$ and $q=r_{M+1}p_{M+1}.$ Let
$P=\left(  m_{n}\right)  ,$ where $m_{n}=$ $\min\operatorname*{supp}x_{n}^{M}$
and $Q=\bigcup\limits_{r_{M+1}\left(  t\right)  =r_{M+1}}\{ m_{n}:x_{n}%
^{M}\subseteq x_{t}^{M+1}\} .$ If $r_{M+1}\left(  t\right)  =r_{M+1},$ let
$G_{t}=\operatorname*{supp}u_{t}^{M+1}.$ Also, let $z_{j}=e_{m_{n_{j}}}$ if
$m_{n_{j}}\in Q.$ Note that if $r_{M+1}\left(  t\right)  =r_{M+1},$
$\theta_{r_{M+1}p_{M+1}}u_{t}^{M+1}$ is an $\mathcal{S}_{q}$-repeated average.
By Lemma \ref{L1pq},
\[
\left\Vert \theta_{r_{M+1}p_{M+1}}u_{t}^{M+1}\right\Vert _{\mathcal{S}%
_{r_{M+1}p_{M+1}-1}}\leq\frac{6}{\min G_{t}}\leq\frac{6}{t}%
\]
if $r_{M+1}\left(  t\right)  =r_{M+1}$ and $\left\Vert \sum\theta
_{r_{M+1}p_{M+1}}u_{t}^{M+1}\right\Vert _{\mathcal{S}_{r_{M+1}p_{M+1}}}\leq6$
for any sum over a finite set of $t$'s such that $r_{M+1}\left(  t\right)
=r_{M+1}.$

Inductively, suppose that the claim is true for some $s<N-M.$ Set $p=p\left(
r_{M+1},...,r_{M+s}\right)  $ and $q=r_{M+s+1}p_{M+s+1}.$ Let $P=\left(
m_{n}\right)  ,$ where $m_{n}=$ $\min\operatorname*{supp}x_{n}^{M+s},$ and
\[
Q=\bigcup\limits_{r_{M+s+1}\left(  t\right)  =r_{M+s+1}}\left\{  m_{n}%
:x_{n}^{M+s}\subseteq x_{t}^{M+s+1},r_{M+s}\left(  n\right)  =r_{M+s}\right\}
.
\]
If $r_{M+s+1}\left(  t\right)  =r_{M+s+1},$ set $G_{t}=\left\{  m_{n}%
:x_{n}^{M+s}\subseteq x_{t}^{M+s+1}\right\}  $. Also$,$\ let $z_{j}%
=\prod_{i=M+1}^{M+s}\theta_{r_{i}p_{i}} u_{n_{j}}^{M+s}\left(  r_{M+1}%
,...,r_{M+s-1}\right)  $ if $m_{n_{j}}\in Q.$ Now%
\begin{align*}
\left\Vert z_{j}\right\Vert _{\ell^{1}}  &  =\left\Vert \prod_{i=M+1}%
^{M+s}\theta_{r_{i}p_{i}}\cdot u_{n_{j}}^{M+s}\left(  r_{M+1},...,r_{M+s-1}%
\right)  \right\Vert _{\ell^{1}}\\
&  =\prod_{i=M+1}^{M+s}\theta_{r_{i}p_{i}}\cdot\left\Vert \left[  x_{n_{j}%
}^{M+s}\left(  r_{M+1},...,r_{M+s-1}\right)  \right]  _{\mathbf{x}^{M}%
}\right\Vert _{\ell^{1}}\leq1
\end{align*}
by Corollary \ref{CorNorm1}. (Note the fact that $L_{i} \geq2$.) By the
inductive hypothesis, $\left\Vert \sum z_{j}\right\Vert _{\mathcal{S}%
_{p\left(  r_{M+1},\dots,r_{M+s}\right)  }}\leq6$ for any sum over a finite
set of $j$'s satisfying $r_{M+s}\left(  n_{j}\right)  =r_{M+s}.$ Finally,
observe that if $r_{M+s+1}\left(  t\right)  =r_{M+s+1}$ and $u_{t}%
^{M+s+1}=\sum_{m_{n}\in G_{t}}c_{n}u_{n}^{M+s},$ then $\theta_{r_{M+s+1}%
p_{M+s+1}}\sum_{m_{n}\in G_{t}}c_{n}e_{m_{n}}$ is an $\mathcal{S}_{q}%
$-repeated average. Thus it follows from Lemma \ref{L1pq} that%

\[
\left\Vert \prod_{i=M+1}^{M+s+1}\theta_{r_{i}p_{i}}\sum u_{t}^{M+s+1}\left(
r_{M+1},\dots,r_{M+s}\right)  \right\Vert _{\mathcal{S}_{p\left(
r_{M+1},\cdots,r_{M+s+1}\right)  }}\leq6
\]
for any sum over a finite set of $t$'s such that $r_{M+s+1}\left(  t\right)
=r_{M+s+1}$ and%
\[
\left\Vert \prod_{i=M+1}^{M+s+1}\theta_{r_{i}p_{i}}u_{t}^{M+s+1}\left(
r_{M+1},\dots,r_{M+s}\right)  \right\Vert _{\mathcal{S}_{p\left(
r_{M+1},\dots,r_{M+s+1}\right)  }}\leq\frac{6}{t}
\]
if $r_{M+s+1}(t) = r_{M+s+1}$. This completes the induction. The Proposition
follows by taking $M+s=N$ and $t=k$ in (\ref{p10}).
\end{proof}

Let $\mathcal{T}$ be an admissible tree and suppose that $0\leq M\leq N-2.$
Say that a collection of nodes $\mathcal{E}$ in $\mathcal{T}$ \textbf{is
subordinated to} $\mathbf{x}^{M}$ if they are pairwise disjoint and for each
$E\in\mathcal{E}$, there exists $j$ such that $E\subseteq x_{j}^{M}.$ Note
that in this case, for every $E\in\mathcal{E}$, there exist unique
$r_{M+1},\dots,r_{N-1}$ such that $E\subseteq x_{k}^{N}\left(  r_{M+1}%
,\dots,r_{N-1}\right)  $. Recall the assumption (\ref{1}) on $k$. Note that if
$x_{j}^{M}\subseteq x_{k}^{N},$ then $j\geq k$ and hence $j$ also satisfies
(\ref{1}) in place of $k.$

\begin{proposition}
\label{PropA}If $\mathcal{E}$ is a collection of nodes in an admissible tree
that is subordinated to $\mathbf{x}^{M}$ and that $o\left(  E\right)
<p\left(  r_{M+1},...,r_{N}\left(  k\right)  \right)  $ for all $E\in
\mathcal{E}$ with $E\subseteq x_{k}^{N}\left(  r_{M+1},...,r_{N-1}\right)  ,$
then
\[
\sum_{E\in\mathcal{E}}t\left(  E\right)  \left\vert \left\vert Ex_{k}%
^{N}\right\vert \right\vert \leq\frac{1}{3N^{2}}.
\]

\end{proposition}

\begin{proof}
Let $\mathcal{E}\left(  r_{M+1},...,r_{N-1}\right)  $ be the set of all nodes
in $\mathcal{E}$ such that $E\subseteq x_{k}^{N}\left(  r_{M+1},...,r_{N-1}%
\right)  $. We have
\begin{align*}
\sum_{E\in\mathcal{E}\left(  r_{M+1},...,r_{N-1}\right)  }  &  t\left(
E\right)  \left\Vert Ex_{k}^{N}\left(  r_{M+1},...,r_{N-1}\right)  \right\Vert
\\
&  \leq\sum_{j\in G}b_{j}\left\Vert x_{j}^{M}\right\Vert \leq\sup_{j\in
G}\left\Vert x_{j}^{M}\right\Vert _{\ell^{1}}\sum_{j\in G}b_{j},
\end{align*}
where $\left(  b_{j}\right)  =\left[  x_{k}^{N}\left(  r_{M+1},...,r_{N-1}%
\right)  \right]  _{\mathbf{x}^{M}}$ and $G$ consists of all $j$'s such that
there exists $E\in\mathcal{E}\left(  r_{M+1},...,r_{N-1}\right)  $ with
$E\subseteq x_{j}^{M}.$ Then
\[
\{ \min\operatorname*{supp}x_{j}^{M}:j\in G\smallsetminus\{ \min G\}\}
\]
is a spreading of a subset of $\{ \min E:E\in\mathcal{E}\} .$ By Lemma
\ref{Tree}, $( x_{j}^{M}) _{j\in G\smallsetminus\left\{  \min G\right\}  }$ is
$\mathcal{S}_{p\left(  r_{M+1},\dots,r_{N}\left(  k\right)  \right)  -1}%
$-admissible. Thus
\[
\sum_{j\in G\smallsetminus\left\{  \min G\right\}  }b_{j}\leq\Vert u_{k}^{N}(
r_{M+1},\dots,r_{N-1}) \Vert_{\mathcal{S}_{p\left(  r_{M+1},\dots,r_{N}\left(
k\right)  \right)  -1}}.
\]
It follows from Proposition \ref{P1} that%
\[
\sum_{j\in G}b_{j}\leq\frac{6}{k}\prod_{i=M+1}^{N}\theta_{r_{i}p_{i}}%
^{-1}+\sup_{j}b_{j}\leq\frac{7}{k}\prod_{i=M+1}^{N}\theta_{r_{i}p_{i}}^{-1}.
\]
Hence, using Corollary \ref{CorNorm},
\begin{align*}
\sum_{E\in\mathcal{E}\left(  r_{M+1},\dots,r_{N-1}\right)  }t  &  \left(
E\right)  \left\Vert Ex_{k}^{N}\left(  r_{M+1},\dots,r_{N-1}\right)
\right\Vert \\
&  \leq\sup_{j\in G}\left\Vert x_{j}^{M}\right\Vert _{\ell^{1}}\frac{7}%
{k}\prod_{i=M+1}^{N}\theta_{r_{i}p_{i}}^{-1}\\
&  \leq2\prod_{i=1}^{M}\theta_{L_{i}p_{i}}^{-1}\cdot\frac{7}{k}\prod
_{i=M+1}^{N}\theta_{r_{i}p_{i}}^{-1}\\
&  \leq\frac{14}{k}\prod_{i=1}^{N}\theta_{L_{i}p_{i}}^{-1}\text{.}%
\end{align*}
Summing over all possible $r_{M+1},\dots,r_{N-1},$ we obtain%
\[
\sum_{E\in\mathcal{E}}t\left(  E\right)  \left\vert \left\vert Ex_{k}%
^{N}\right\vert \right\vert \leq\frac{14}{k}\prod_{i=1}^{N}L_{i}\theta
_{L_{i}p_{i}}^{-1}\leq\frac{1}{3N^{2}}%
\]
\newline by (\ref{1}).
\end{proof}

Next, consider a set of nodes $\mathcal{E}^{\prime}$ in $\mathcal{T}$ that is
subordinated to $\mathbf{x}^{M}$ and that
\[
o\left(  E\right)  \geq p\left(  r_{M+1}+1,r_{M+2},\dots,r_{N}\left(
k\right)  \right)
\]
for all $E\in\mathcal{E}^{\prime}$ with $E\subseteq x_{k}^{N}\left(
r_{M+1},\dots,r_{N-1}\right)  .$ In analogy to the above, for given
$r_{M+1},\dots,r_{N-1},$ let $\mathcal{E}^{\prime}\left(  r_{M+1}%
,\dots,r_{N-1}\right)  $ be the set of all nodes in $\mathcal{E}^{\prime}$
such that $E\subseteq x_{k}^{N}\left(  r_{M+1},\dots,r_{N-1}\right)  $.

\begin{proposition}
\label{12}$\sum_{E\in\mathcal{E}^{\prime}}t( E) \Vert Ex_{k}^{N}\Vert\leq
\frac{1}{3N^{2}}\Vert x_{k}^{N}\Vert_{X_{M}}.$
\end{proposition}

\begin{proof}
We have%
\begin{align*}
\sum_{E\in\mathcal{E}^{\prime}(r_{M+1},\dots,r_{N-1})}  &  t(E)\Vert
Ex_{k}^{N}(r_{M+1},\dots,r_{N-1})\Vert\\
&  \leq\theta_{p\left(  r_{M+1}+1,\dots,r_{N}\left(  k\right)  \right)  }%
\sum_{j\in G^{\prime}}a_{j}\Vert x_{j}^{M}\Vert\\
&  \leq\theta_{p\left(  r_{M+1}+1,\dots,r_{N}\left(  k\right)  \right)  }%
\sup_{j\in G^{\prime}}\Vert x_{j}^{M}\Vert_{\ell^{1}}\sum_{j\in G^{\prime}%
}a_{j}%
\end{align*}
where $(a_{j})=[x_{k}^{N}(r_{M+1},\dots,r_{N-1})]_{\mathbf{x}^{M}}$ and
$G^{\prime}$ consists of all $j$'s such that there exists $E\in\mathcal{E}%
^{\prime}\left(  r_{M+1},\dots,r_{N-1}\right)  $ with $E\subseteq x_{j}^{M}.$
But
\begin{align*}
\sum_{j\in G^{\prime}}a_{j}  &  \leq\Vert\lbrack x_{k}^{N}(r_{M+1}%
,\dots,r_{N-1})]_{\mathbf{x}^{M}}\Vert_{\ell^{1}}\\
&  \leq2\theta_{r_{N}(k)p_{N}}^{-1}\prod_{i=M+1}^{N-1}\theta_{r_{i}p_{i}}%
^{-1}L_{i}^{-1}\quad\text{by Corollary \ref{CorNorm1}.}%
\end{align*}
Applying Corollary \ref{CorNorm} to the above, we have%
\begin{align}
\sum_{E\in\mathcal{E}^{\prime}\left(  r_{M+1},\dots,r_{N-1}\right)  }  &
t\left(  E\right)  \Vert Ex_{k}^{N}\left(  r_{M+1},\dots,r_{N-1}\right)
\Vert\nonumber\\
&  \leq4\theta_{p\left(  r_{M+1}+1,\dots,r_{N}\left(  k\right)  \right)
}\theta_{r_{N}\left(  k\right)  p_{N}}^{-1}\prod_{i=1}^{M}\theta_{L_{i}p_{i}%
}^{-1}\prod_{i=M+1}^{N-1}\theta_{r_{i}p_{i}}^{-1}L_{i}^{-1}. \label{B}%
\end{align}
Recall the lower estimate for $\Vert x_{k}^{N}\Vert_{X_{M}}$ given by
(\ref{B'}) in Corollary \ref{ModNormM}. For fixed $r_{M+1},\dots,r_{N-1},$ the
ratio of $\left(  \ref{B}\right)  $ with the $\left(  r_{M+1},\dots
,r_{N-1}\right)  $-indexed term in $(\ref{B'})$ is
\begin{align*}
&  \leq\frac{8}{\theta_{1}}\frac{\theta_{p\left(  r_{M+1}+1,\dots,r_{N}\left(
k\right)  \right)  }}{\theta_{p\left(  r_{M+1},\dots,r_{N}\left(  k\right)
\right)  }}\prod_{i=1}^{M}\theta_{L_{i}p_{i}}^{-1}\\
&  =\frac{8}{\theta_{1}}\frac{\theta_{p_{M+1}+p\left(  r_{M+1},\dots
,r_{N}\left(  k\right)  \right)  }}{\theta_{p\left(  r_{M+1},\dots
,,r_{N}\left(  k\right)  \right)  }}\prod_{i=1}^{M}\theta_{L_{i}p_{i}}^{-1}\\
&  \leq\frac{8}{\theta_{1}}\frac{\theta_{1}}{24N^{2}}\prod_{i=1}^{M}%
\theta_{L_{i}p_{i}}\prod_{i=1}^{M}\theta_{L_{i}p_{i}}^{-1}\quad\text{by
Condition \textbf{(A)},}\\
&  =\frac{1}{3N^{2}}.
\end{align*}
Hence
\[
\sum_{E\in\mathcal{E}^{\prime}}t(E)\Vert Ex_{k}^{N}\Vert\leq\frac{1}{3N^{2}%
}\Vert x_{k}^{N}\Vert_{X_{M}}.
\]

\end{proof}

In the next two results, let $( d_{j}) =[ x_{k}^{N}] _{\mathbf{x}^{M+1}}.$
Recall the convention that $x^{N}_{k}( r_{M+2},\dots,r_{N-1}) = x^{N}_{k}$ if
$M = N-2$.

\begin{lemma}
\label{Lcolorskipped} Suppose that $0 \leq M \leq N-2$. Given $r_{M+2}%
,\dots,r_{N-1}$, write
\[
K=\{ j:x_{j}^{M+1}\subseteq x^{N}_{k}( r_{M+2},\dots,r_{N-1})\} .
\]
If $J$ is an $L_{M+1}$-skipped set, then
\begin{equation}
\sum_{j\in J\cap K}d_{j}\leq\left(  L_{M+1}^{-1}+k^{-1}\right)  \sum_{j\in
K}d_{j}\leq\frac{3}{2}L_{M+1}^{-1}\sum_{j\in K}d_{j}. \label{LcolorskippedCor}%
\end{equation}

\end{lemma}

\begin{proof}
The second inequality follows from the choice of $k$ since $k\geq2L_{M+1}$ by
(\ref{1}). Recall the notation from $\boldsymbol{(\beta)}$ expressing
\[
x_{i}^{M+2}=\sum_{j\in I_{i}^{M+2}}a_{j}x_{j}^{M+1}.
\]
For each $i$ such that $x_{i}^{M+2}\subseteq x_{k}^{N},$ let $J_{i}=J\cap
I_{i}^{M+2}.$ Then $J_{i}$ is an $L_{M+1}$-skipped subset of the integer
interval $I_{i}^{M+2}.$ By Lemma \ref{Lskippedseq},%
\begin{align*}
\sum_{j\in J_{i}}d_{j}  &  \leq L_{M+1}^{-1}\sum_{j\in I_{i}^{M+2}}d_{j}%
+\sup_{j\in I_{i}^{M+2}}d_{j}\\
&  =L_{M+1}^{-1}\sum_{j\in I_{i}^{M+2}}d_{j}+d_{\min I_{i}^{M+2}}.
\end{align*}
Now $J\cap K=\cup_{i\in K^{\prime}}J_{i},$ where $K^{\prime}=\{i:x_{i}%
^{M+2}\subseteq x_{k}^{N}(r_{M+2},\dots,r_{N-1})\}.$ Thus
\begin{align*}
\sum_{j\in J\cap K}d_{j}  &  \leq{L_{M+1}^{-1}}\sum_{i\in K^{\prime}}%
\sum_{j\in I_{i}^{M+2}}d_{j}+\sum_{i\in K^{\prime}}d_{\min I_{i}^{M+2}}\\
&  ={L_{M+1}^{-1}}\sum_{j\in K}d_{j}+\sum_{i\in K^{\prime}}d_{\min I_{i}%
^{M+2}}.
\end{align*}
If $[x_{k}^{N}]_{\mathbf{x}^{M+2}}=\left(  b_{i}\right)  ,$ then for all $j\in
I_{i}^{M+2},$ we can express $d_{j}=b_{i}a_{j},$ where $\theta_{r_{M+2}\left(
i\right)  p_{M+2}}\sum_{j\in I_{i}^{M+2}}a_{j}e_{m_{j}}$ is an \thinspace
$\mathcal{S}_{r_{M+2}\left(  i\right)  p_{M+2}}$-repeated average, with
$e_{m_{j}}=\min\operatorname*{supp}x_{j}^{M+1}.$ In particular,
\begin{align*}
\theta_{r_{M+2}\left(  i\right)  p_{M+2}}a_{j_{0}}  &  \leq i^{-1}\leq
k^{-1}\\
&  =k^{-1}\theta_{r_{M+2}\left(  i\right)  p_{M+2}}\sum_{j\in I_{i}^{M+2}%
}a_{j}%
\end{align*}
for all $j_{0}\in I_{i}^{M+2}.$ Thus
\begin{align*}
d_{\min I_{i}^{M+2}}  &  =b_{i}a_{\min I_{i}^{M+2}}\leq b_{i}k^{-1}\sum_{j\in
I_{i}^{M+2}}a_{j}\\
&  \leq k^{-1}\sum_{j\in I_{i}^{M+2}}b_{i}a_{j}=k^{-1}\sum_{j\in I_{i}^{M+2}%
}d_{j}.
\end{align*}
Therefore,%
\[
\sum_{i\in K^{\prime}}d_{\min I_{i}^{M+2}}\leq k^{-1}\sum_{i\in K^{\prime}%
}\sum_{j\in I_{i}^{M+2}}d_{j}=k^{-1}\sum_{j\in K}d_{j}.
\]
Hence
\[
\sum_{j\in J\cap K}d_{j}=(L_{M+1}^{-1}+k^{-1})\sum_{j\in K}d_{j}.
\]

\end{proof}

We say that an admissible tree $\mathcal{T}$ \textbf{is subordinated to}
$\mathbf{x}^{M}$ if its set of base nodes is subordinated to $\mathbf{x}^{M}$
and any leaf that is not at the base is a singleton. Given an admissible tree
that is subordinated to $\mathbf{x}^{M},$ let $\mathcal{E}^{\prime\prime}$ be
the collection of all base nodes $E$ in $\mathcal{T}$ such that $p\left(
r_{M+1},...,r_{N}\left(  k\right)  \right)  \leq o\left(  E\right)  <p\left(
r_{M+1}+1,...,r_{N}\left(  k\right)  \right)  $ if $E\subseteq x_{k}%
^{N}\left(  r_{M+1},...,r_{N-1}\right)  .$ It follows from Condition
\textbf{(B)} that for $E\in\mathcal{E}^{\prime\prime}$, $o\left(  E\right)  $
uniquely determines $r_{M+1},...,r_{N-1}$ such that $E\subseteq x_{k}%
^{N}(r_{M+1},...,r_{N-1}).$ Let $\mathcal{D}$ denote the set of all $D$'s that
are immediate predecessors of some $E\in\mathcal{E}^{\prime\prime}.$ We say
that $D$ \textbf{effectively intersects} $x_{j}^{M+1}$ for some $j$ if there
exists $E\in\mathcal{E}^{\prime\prime}$ such that $E\subseteq
\operatorname*{supp}x_{j}^{M+1}\cap D$. Let $\mathcal{\tilde{D}}$ be the
subcollection of all $D\in\mathcal{D}$ such that $D$ effectively intersects at
least two $x_{j}^{M+1}$'s. For each $D\in\mathcal{\tilde{D}}$, let
$J(D)=\{j:D\text{ effectively intersects }x_{j}^{M+1}\}$, then $J\left(
D\right)  $ is an $L_{M+1}$-skipped set. Indeed, if $D\in\mathcal{D}$ and
$E_{1},E_{2}$ are successors of $D$ in $\mathcal{E}^{\prime\prime}$ such that
$E_{i}\subseteq\operatorname*{supp}x_{j_{i}}^{M+1}\cap D$, $i=1,2,$ and
$j_{1}<j_{2},$ then $o\left(  E_{1}\right)  =o\left(  E_{2}\right)  $ and
hence $r_{M+1}\left(  j_{1}\right)  =r_{M+1}\left(  j_{2}\right)  .$ Thus
$j_{2}-j_{1}\geq L_{M+1}.$ Let $J=\cup_{D\in\mathcal{\tilde{D}}}J\left(
D\right)  .$ If the elements of $\mathcal{\tilde{D}}$ are arranged in order,
then the union of $J\left(  D\right)  $ taken over every other $D\in
\mathcal{\tilde{D}}$ is an $L_{M+1}$-skipped set. Hence $J$ is the union of at
most two $L_{M+1}$-skipped sets.

\begin{proposition}
\label{14}
\[
\sum_{D\in\mathcal{\tilde{D}}}\sum_{\substack{E\in\mathcal{E}^{^{\prime\prime
}}\\\!E\subseteq D}}t\left(  E\right)  \Vert Ex_{k}^{N}\Vert\leq\frac
{1}{3N^{2}}\Vert x_{k}^{N}\Vert_{X_{M}}.
\]

\end{proposition}

\begin{proof}
Let $(d_{j})$ be as in Lemma \ref{Lcolorskipped} and $g(j)=p(r_{M+1}%
,\dots,r_{N}(k))$ if $x_{j}^{M+1}\subseteq x_{k}^{N}(r_{M+1},\dots,r_{N-1}).$
Then
\[
\sum_{D\in\mathcal{\tilde{D}}}\sum_{\substack{E\in\mathcal{E}^{^{\prime\prime
}}\\\!E\subseteq D}}t(E)\Vert Ex_{k}^{N}\Vert\leq\sum_{D\in\mathcal{\tilde{D}%
}}\sum_{j\in J(D)}\theta_{g(j)}d_{j}\Vert x_{j}^{M+1}\Vert_{\ell^{1}}.
\]
But $x_{j}^{M+1}=\sum_{\ell\in I_{j}^{M+1}}a_{\ell}x_{\ell}^{M}$ with $\sum
a_{\ell}=\theta_{r_{M+1}\left(  j\right)  p_{M+1}}^{-1}.$ Hence
\begin{align}
\sum_{D\in\mathcal{\tilde{D}}}\sum_{\substack{E\in\mathcal{E}^{^{\prime\prime
}}\\\!E\subseteq D}}t\left(  E\right)  \Vert Ex_{k}^{N}\Vert &  \leq\sup
_{\ell}\Vert x_{\ell}^{M}\Vert_{\ell^{1}}\sum_{D\in\mathcal{\tilde{D}}}%
\sum_{j\in J\left(  D\right)  }\theta_{g\left(  j\right)  }d_{j}%
\theta_{r_{M+1}\left(  j\right)  p_{M+1}}^{-1}\nonumber\\
&  \leq2\sup_{\ell}\Vert x_{\ell}^{M}\Vert_{\ell^{1}}\sum_{j\in J}%
\theta_{g\left(  j\right)  }d_{j}\theta_{r_{M+1}\left(  j\right)  p_{M+1}%
}^{-1}, \label{star}%
\end{align}
since each $j$ belongs to at most two $J\left(  D\right)  .$ Fix
$r_{M+2},...,r_{N-1}$ and let $K$ be as in Lemma \ref{Lcolorskipped}. Since
$J$ is the union of at most two $L_{M+1}$-skipped sets,
\begin{align*}
\sum_{j\in J\cap K}\theta_{g\left(  j\right)  }d_{j}\theta_{r_{M+1}\left(
j\right)  p_{M+1}}^{-1}  &  \leq\sup_{j\in K}\frac{\theta_{g\left(  j\right)
}}{\theta_{r_{M+1}\left(  j\right)  p_{M+1}}}\sum_{j\in J\cap K}d_{j}\\
&  \leq\frac{3}{L_{M+1}}\sup_{j\in K}\frac{\theta_{g\left(  j\right)  }%
}{\theta_{r_{M+1}\left(  j\right)  p_{M+1}}}\sum_{j\in K}d_{j}\quad\text{by
(\ref{LcolorskippedCor}).}%
\end{align*}
However,
\begin{align*}
\sup_{j\in K}\frac{\theta_{g\left(  j\right)  }}{\theta_{r_{M+1}\left(
j\right)  p_{M+1}}}  &  \leq\sup_{1\leq j\leq L_{M+1}}\frac{\theta
_{r_{M+1}\left(  j\right)  p_{M+1}+p\left(  r_{M+2},...,r_{N}\left(  k\right)
\right)  }}{\theta_{r_{M+1}\left(  j\right)  p_{M+1}}}\\
&  \leq F\left(  L_{M+1}\right)  \sum_{r_{M+1}}\frac{\theta_{r_{M+1}%
p_{M+1}+p\left(  r_{M+2},...,r_{N}\left(  k\right)  \right)  }}{\theta
_{r_{M+1}p_{M+1}}}.
\end{align*}
by Condition $(\ddag)$. Therefore,
\[
\sum_{j\in J\cap K}\theta_{g\left(  j\right)  }d_{j}\theta_{r_{M+1}\left(
j\right)  p_{M+1}}^{-1}\leq\frac{3F\left(  L_{M+1}\right)  }{L_{M+1}}%
\sum_{r_{M+1}}\frac{\theta_{p\left(  r_{M+1},\dots,r_{N}\left(  k\right)
\right)  }}{\theta_{r_{M+1}p_{M+1}}}\sum_{j\in K}d_{j}.
\]
Note that
\[
\sum_{j\in K}d_{j}=\Vert\left[  x_{k}^{N}\left(  r_{M+2},\dots,r_{N-1}\right)
\right]  _{\mathbf{x}^{M+1}}\Vert_{\ell^{1}}\leq2\theta_{r_{N}\left(
k\right)  p_{N}}^{-1}\prod_{i=M+2}^{N-1}\theta_{r_{i}p_{i}}^{-1}L_{i}^{-1}%
\]
by Corollary \ref{CorNorm1}. Summing over all $r_{M+2},\dots,r_{N-1},$ we have%
\begin{align}
&  \sum_{j\in J}\theta_{g\left(  j\right)  }d_{j}\theta_{r_{M+1}\left(
j\right)  p_{M+1}}^{-1}\nonumber\\
&  \leq\frac{6F\left(  L_{M+1}\right)  }{L_{M+1}}\sum_{r_{M+1},\dots,r_{N-1}%
}\frac{\theta_{p\left(  r_{M+1},\dots,r_{N-1},r_{N}\left(  k\right)  \right)
}}{\theta_{r_{M+1}p_{M+1}}}\cdot\theta_{r_{N}\left(  k\right)  p_{N}}%
^{-1}\prod_{i=M+2}^{N-1}\theta_{r_{i}p_{i}}^{-1}L_{i}^{-1}\nonumber\\
&  =6F\left(  L_{M+1}\right)  \sum_{r_{M+1},\dots,r_{N-1}}\theta_{p\left(
r_{M+1},\dots,r_{N-1},r_{N}\left(  k\right)  \right)  }\theta_{r_{N}\left(
k\right)  p_{N}}^{-1}\prod_{i=M+1}^{N-1}\theta_{r_{i}p_{i}}^{-1}\text{ }%
L_{i}^{-1}. \label{E10}%
\end{align}
Comparing (\ref{star}) and (\ref{E10}) with (\ref{B'}) in Corollary
\ref{ModNormM}, we see that
\begin{align*}
\sum_{D\in\mathcal{\tilde{D}}}\sum_{\substack{E\in\mathcal{E}^{^{\prime\prime
}}\\\!E\subseteq D}}t(E)\Vert Ex\Vert &  \leq24\theta_{1}^{-1}F(L_{M+1})\Vert
x_{k}^{N}\Vert_{X_{M}}\sup_{\ell}\Vert x_{\ell}^{M}\Vert_{\ell^{1}}\\
&  \leq48F\left(  L_{M+1}\right)  \theta_{1}^{-1}\prod_{i=1}^{M}\theta
_{L_{i}p_{i}}^{-1}\Vert x_{k}^{N}\Vert_{X_{M}}\\
&  \quad\quad\text{ by Corollary \ref{CorNorm},}\\
&  \leq\frac{1}{3N^{2}}\Vert x_{k}^{N}\Vert_{X_{M}}%
\end{align*}
by condition \textbf{(C)}.
\end{proof}

\begin{definition}
Given $N,p\in\mathbb{N}$ define
\[
\Theta_{p} = \Theta_{p}\left(  N\right)  =\max\bigl\{
\prod_{i=1}^{N}\theta_{\ell_{i}}: \ell_{i}\in{\mathbb{N}}, \sum_{i=1}^{N}%
\ell_{i}=p\bigr\}  .
\]

\end{definition}

For any $N\in\mathbb{N}$ and $V\in[ \mathbb{N}] $, choose integer sequences
$\left(  p_{k}\right)  _{k=1}^{N}$ and $\left(  L_{k}\right)  _{k=1}^{N},$ and
sequences of vectors $\mathbf{x}^{0},\mathbf{x}^{1},\dots,\mathbf{x}^{N}$ as above.

\begin{theorem}
\label{Th21}There exists a finitely supported vector $x\in\operatorname*{span}%
\left\{  e_{k}:k\in V\right\}  $ such that%
\begin{equation}
\left\Vert x\right\Vert \leq\bigl(\frac{2}{N}+4\theta_{1}^{-1}\sup
_{r_{1},\dots,r_{N-1}}\frac{\Theta_{p\left(  r_{1},\dots,r_{N}\left(
k\right)  \right)  }}{\theta_{p\left(  r_{1},\dots,r_{N}\left(  k\right)
\right)  }}\bigr)\left\Vert x\right\Vert _{X_{M}}. \label{Th21E}%
\end{equation}

\end{theorem}

\begin{proof}
Consider an admissible tree $\mathcal{T}$ that is subordinated to
$\mathbf{x}^{M}$, $0\leq M\leq N-2.$ Let $\mathcal{E}$ and $\mathcal{E}%
^{\prime}$ be the set of all base nodes such that $o\left(  E\right)
<p\left(  r_{M+1},\dots,r_{N}\left(  k\right)  \right)  $, respectively,
$o\left(  E\right)  \geq p\left(  r_{M+1}+1,\dots,r_{K}\left(  k\right)
\right)  $ if $E\subseteq x_{k}^{M}\left(  r_{M+1},\dots,r_{N-1}\right)  .$
Also, define $\mathcal{E}^{\prime\prime},\mathcal{D}$ and $\mathcal{\tilde{D}%
}$ as in the discussion preceding Proposition \ref{14}. Finally, let
$\mathcal{E}^{\prime\prime\prime}$ be the set of all leaves of $\mathcal{T}$
not at the base. By Proposition \ref{14},
\begin{align*}
\sum_{E\in\mathcal{E}^{^{\prime\prime}}}t\left(  E\right)  \Vert Ex_{k}%
^{N}\Vert &  = \sum_{D\in\mathcal{D}}\sum_{\substack{E\in\mathcal{E}%
^{^{\prime\prime}}\\\!E\subseteq D}}t\left(  E\right)  \Vert Ex_{k}^{N}\Vert\\
&  \leq\frac{1}{3N^{2}}\Vert x_{k}^{N}\Vert_{X_{M}}+\sum_{D\in
\mathcal{D\smallsetminus\tilde{D}}}\sum_{\substack{E\in\mathcal{E}%
^{^{\prime\prime}}\\\!E\subseteq D}}t\left(  E\right)  \Vert Ex_{k}^{N}\Vert.
\end{align*}
If $D\in\mathcal{D\smallsetminus\tilde{D}}$, $D$ effectively intersects at
most one $x_{j}^{M+1}.$ Set $D^{\prime}=D\cap\operatorname*{supp}x_{j}^{M+1}$
($D^{\prime}=\emptyset$ if no such $j$ exists). Then
\[
\sum_{\substack{E\in\mathcal{E}^{^{\prime\prime}}\\\!E\subseteq D}}t\left(
E\right)  \Vert Ex_{k}^{N}\Vert=\sum_{\substack{E\in\mathcal{E}^{^{\prime
\prime}}\\E\subseteq D^{\prime}}}t\left(  E\right)  \Vert Ex_{k}^{N}\Vert\leq
t\left(  D\right)  \Vert D^{\prime}x_{k}^{N}\Vert.
\]
Now let $\mathcal{T}^{\prime}$ be a tree obtained from $\mathcal{T}$ by taking
all $D\in\mathcal{D\smallsetminus}\mathcal{\tilde{D}}$, all $E\in
\mathcal{E}^{\prime\prime\prime}$ and all their ancestors, with each
$D\in\mathcal{D\smallsetminus}\mathcal{\tilde{D}}$ modified into $D^{\prime}$
as described above. Then $\mathcal{T}^{\prime}$ is an admissible tree that is
subordinated to $\mathbf{x}^{M+1}$ and $H\left(  \mathcal{T}^{\prime}\right)
<H\left(  \mathcal{T}\right)  .$ (Note that every node in $\mathcal{E}%
^{\prime\prime\prime}$ is a singleton.) By Propositions \ref{PropA} and
\ref{12} and the above,%
\begin{align*}
\sum_{E \in\mathcal{L}(\mathcal{T})}t(E)\|Ex_{k}^{N}\|  &  \leq\frac{1}%
{3N^{2}}+\frac{1}{3N^{2}}\Vert x_{k}^{N}\Vert_{X_{M}}+\left(  \sum
_{E\in\mathcal{E}^{\prime\prime}}+\sum_{E\in\mathcal{E}^{\prime\prime\prime}%
}\right)  t\left(  E\right)  \Vert Ex_{k}^{N}\Vert\\
&  \leq\frac{1}{N^{2}}{\Vert x_{k}^{N}\Vert_{X_{M}}}+\sum_{D\in
\mathcal{D\smallsetminus\tilde{D}}}t\left(  D\right)  \Vert D^{\prime}%
x_{k}^{N}\Vert+\sum_{E\in\mathcal{E}^{\prime\prime\prime}}t\left(  E\right)
\Vert Ex_{k}^{N}\Vert\\
&  =\frac{1}{N^{2}}{\Vert x_{k}^{N}\Vert_{X_{M}}}+\sum_{E \in\mathcal{L}%
(\mathcal{T^{\prime}})}t(E)\|Ex_{k}^{N}\| .
\end{align*}

Now let $\mathcal{T}$ be an admissible tree all of whose leaves are
singletons. Let $\mathcal{T}_{1}$ be the subtree of $\mathcal{T}$ consisting
of leaves $E$ in $\mathcal{T}$ with $h(E)<N$ and their ancestors. Then
$\mathcal{T}_{1}$ is subordinated to $\mathbf{x}^{0}$ and $H\left(
\mathcal{T}_{1}\right)  \leq N-1.$ By the above argument, there is an
admissible tree $\mathcal{T}_{1}^{\prime}$ respecting $\mathbf{x}^{1}$ with
${H}(\mathcal{T}_{1}^{\prime})\leq N-2$ so that
\[
\sum_{E\in\mathcal{L}(\mathcal{T}_{1})}t(E)\Vert Ex_{k}^{N}\Vert\leq\sum
_{E\in\mathcal{L}(\mathcal{T}_{1}^{\prime})}t(E)\Vert Ex_{k}^{N}\Vert+\frac
{1}{N^{2}}{\Vert x_{k}^{N}\Vert_{X_{M}}}.
\]
Repeating the argument, we reach an admissible tree $\mathcal{T}_{1}^{(N-1)}$
subordinated to $\mathbf{x}^{N-1}$ with $H(\mathcal{T}_{1}^{(N-1)})=0$ such
that
\[
\sum_{E\in\mathcal{L}(\mathcal{T}_{1})}t(E)\Vert Ex_{k}^{N}\Vert\leq\sum
_{E\in\mathcal{L}(\mathcal{T}_{1}^{(N-1)})}t(E)\Vert Ex_{k}^{N}\Vert
+\frac{N-1}{N^{2}}{\Vert x_{k}^{N}\Vert_{X_{M}}}.
\]
Since $H(\mathcal{T}_{1}^{(N-1)})=0$ and $\mathcal{T}_{1}^{(N-1)}$ is
subordinated to $\mathbf{x}^{N-1},$ $\mathcal{T}_{1}^{(N-1)}$ consists of a
single node $E$ such that $E\subseteq x_{j_{0}}^{N-1}$ for some $j_{0}.$
Recall that $x_{k}^{N}=\sum_{j\in I_{k}^{N}}a_{j}x_{j}^{N-1},$ where
$0\leq\theta_{r_{N}\left(  k\right)  p_{N}}a_{j}\leq k^{-1}$ for all $j\in
I_{k}^{N}.$ Hence%
\begin{align*}
\sum_{E\in\mathcal{L}(\mathcal{T}_{1}^{(N-1)})}t(E)\Vert Ex_{k}^{N}\Vert &
\leq a_{j_{0}}\Vert x_{j_{0}}^{N-1}\Vert_{\ell^{1}}\\
&  \leq2\theta_{r_{N}\left(  k\right)  p_{N}}^{-1}k^{-1}\prod_{i=1}%
^{N-1}\theta_{L_{i}p_{i}}^{-1}\quad\text{ by Corollary \ref{CorNorm},}\\
&  \leq\frac{1}{N^{2}}\quad\text{ by (\ref{1}).}%
\end{align*}
Therefore,
\begin{equation}
\sum_{E\in\mathcal{L}(\mathcal{T}_{1})}t(E)\Vert Ex_{k}^{N}\Vert\leq\frac
{1}{N}\Vert x_{k}^{N}\Vert_{X_{M}}. \label{i}%
\end{equation}
Let $\mathcal{T}_{2}$ be the subtree of $\mathcal{T}$ consisting of leaves $E$
in $\mathcal{T}$ with $h(E)\geq N$ and their ancestors. Since every leaf in
$\mathcal{T}_{2}$ is a singleton, the set of all leaves is subordinated to
$\mathbf{x}^{0}.$ Let $\mathcal{G}$ be the collection of all leaves $E$ of
$\mathcal{T}_{2}$ such that $o\left(  E\right)  <p\left(  r_{1},\dots
,r_{N}\left(  k\right)  \right)  $ if $E\subseteq x_{k}^{N}\left(  r_{1}%
,\dots,r_{N}\left(  k\right)  \right)  .$ Then
\[
\sum_{E\in\mathcal{G}}t\left(  E\right)  \Vert Ex_{k}^{N}\Vert\leq\frac
{1}{3N^{2}}\quad\text{ by Proposition \ref{PropA}.}%
\]
Hence
\[
\sum_{E\in\mathcal{L}(\mathcal{T}_{2})}t(E)\Vert Ex_{k}^{N}\Vert\leq\frac
{1}{3N^{2}}+\sum_{E\in\mathcal{G}^{\prime}}t\left(  E\right)  \Vert Ex_{k}%
^{N}\Vert,
\]
where $\mathcal{G}^{\prime}$ consists of all leaves of $\mathcal{T}_{2}$ that
are not in $\mathcal{G}$. If $E\in\mathcal{G}^{\prime}$ and $E\subseteq
x_{k}^{N}\left(  r_{1},\dots,r_{N-1}\right)  ,$ then $o\left(  E\right)  \geq
p\left(  r_{1},\dots,r_{N}\left(  k\right)  \right)  $ and $h\left(  E\right)
\geq N.$ Thus $t\left(  E\right)  =\prod_{i=1}^{j}\theta_{\ell_{j}}$ with
$j\geq N$ and $\sum\ell_{j}\geq p\left(  r_{1},\dots,r_{N}\left(  k\right)
\right)  .$ Since $\left(  \theta_{n}\right)  $ is regular and decreasing,
$t\left(  E\right)  \leq\Theta_{p\left(  r_{1},\dots,r_{N}\left(  k\right)
\right)  }.$ Therefore, using the estimates from Corollary \ref{CorNorm1} and
Proposition \ref{ModNorm}, we have
\begin{align}
\sum_{E\in\mathcal{L}(\mathcal{T}_{2})}t(E)\Vert Ex_{k}^{N}\Vert &  \leq
\frac{1}{3N^{2}}+\sum_{r_{1},\dots,r_{N-1}}\sum_{\substack{E\in\mathcal{G}%
^{\prime}\\E\subseteq x_{k}^{N}\left(  r_{1},\dots,r_{N-1}\right)  }}t\left(
E\right)  \Vert Ex_{k}^{N}\Vert\nonumber\\
&  \leq\frac{1}{3N^{2}}+\sum_{r_{1},\dots,r_{N-1}}\Theta_{p\left(  r_{1}%
,\dots,r_{N}\left(  k\right)  \right)  }\Vert x_{k}^{N}\left(  r_{1}%
,\dots,r_{N-1}\right)  \Vert_{\ell^{1}}\nonumber\\
&  \leq\bigl(\frac{1}{3N^{2}}+4\theta_{1}^{-1}\sup_{r_{1},\dots,r_{N-1}}%
\frac{\Theta_{p\left(  r_{1},\dots,r_{N}\left(  k\right)  \right)  }}%
{\theta_{p\left(  r_{1},\dots,r_{N}\left(  k\right)  \right)  }}\bigr)\Vert
x_{k}^{N}\Vert_{X_{M}} \label{ii}%
\end{align}
Combining (\ref{i}) and (\ref{ii}) and maximizing over all admissible trees
gives%
\begin{align*}
\Vert x_{k}^{N}\Vert &  =\max_{\mathcal{T}}\mathcal{T}x_{k}^{N}\\
&  \leq\bigl(\frac{2}{N}+4\theta_{1}^{-1}\sup_{r_{1},...,r_{N-1}}\frac
{\Theta_{p\left(  r_{1},...,r_{N}\left(  k\right)  \right)  }}{\theta
_{p\left(  r_{1},...,r_{N}\left(  k\right)  \right)  }}\bigr)\Vert x_{k}%
^{N}\Vert_{X_{M}}.
\end{align*}

\end{proof}

\section{\noindent Main Results}

Recall that we define $\theta= \lim\theta_{n}^{1/n} =\sup\theta_{n}^{1/n}$ for
a regular sequence $(\theta_{n})$. Also set $\varphi_{n} = \theta_{n}%
/\theta^{n}$. It was mentioned in the discussion at the beginning of \S 2 that
$X$ and $X_{M}$ are not isomorphic if $\theta= 1$. If $\theta< 1$ and
$\varphi_{N} = 1$ for some $N$, then $X$ and $X_{M}$ are isomorphic by
Proposition \ref{Prop5}. We shall presently show that $X$ and $X_{M}$ are not
isomorphic under some mild conditions on $(\varphi_{n})$. For the remainder of
the section, assume that $\theta< 1$.

\begin{proposition}
\noindent\label{III}If $\inf\varphi_{n}=c>0.$ Then $\left(  \theta_{n}\right)
$ satisfies $\left(  \lnot\dag\right)  $ and $\left(  \ddag\right)  $.
\end{proposition}

\begin{proof}
Indeed,%
\[
\frac{\theta_{m+n}}{\theta_{n}}=\frac{\varphi_{m+n}}{\varphi_{n}}\theta
^{m}\leq\frac{1}{c}\theta^{m}\text{ for all }m,n\in\mathbb{N}\text{.}%
\]
Thus $\left(  \lnot\dag\right)  $ holds. Also,%
\[
\sum_{i=1}^{R}\frac{\theta_{s_{i}+t}}{\theta_{s_{i}}}=\sum_{i=1}^{R}%
\frac{\varphi_{s_{i}+t}}{\varphi_{s_{i}}}\theta^{t}\geq cR\theta^{t}.
\]
On the other hand,
\[
\max_{1\leq i\leq R}\frac{\theta_{s_{i}+t}}{\theta_{s_{i}}}=\max_{1\leq i\leq
R}\frac{\varphi_{s_{i}+t}}{\varphi_{s_{i}}}\theta^{t}\leq\frac{\theta^{t}}%
{c}.
\]
Thus $\left(  \ddag\right)  $ holds with $F\left(  R\right)  =\frac{1}{c^{2}%
R}.$
\end{proof}

\begin{theorem}
\label{Th16}\noindent If $0<c=\inf\varphi_{n}\leq\sup\varphi_{n}=d<1,$ then
$X$ is not isomorphic to $X_{M}.$
\end{theorem}

\begin{proof}
Let $\varepsilon>0$ and $V\in[\mathbb{N}]$ be given. Choose $N\in\mathbb{N}$
such that $\frac{2}{N}+4\theta_{1}^{-1}\frac{d^{N}}{c}<\varepsilon.$ Obtain
from Theorem \ref{Th21} a vector $x\in\operatorname*{span}\left\{  e_{k}:k\in
V\right\}  $ that satisfies (\ref{Th21E}). Let $p\in\mathbb{N}$, if $\left(
\ell_{i}\right)  _{i=1}^{N}$ is a sequence of positive integers such that
$\sum_{i=1}^{N}\ell_{i}=p,$ then
\[
\prod_{i=1}^{N}\theta_{\ell_{i}}=\theta^{p}\prod_{i=1}^{N}\varphi_{\ell_{i}%
}\leq\theta^{p}d^{N}%
\]
and
\[
\theta_{p}=\varphi_{p}\theta^{p}\geq c\theta^{p}.
\]
Thus
\[
\sup_{p}\frac{\Theta_{p}}{\theta_{p}}\leq\frac{d^{N}}{c}.
\]
It follows from (\ref{Th21E}) that
\[
\Vert x\Vert\leq\bigl(\frac{2}{N}+4\theta_{1}^{-1}\frac{d^{N}}{c}\bigr)\Vert
x\Vert_{X_{M}}<\varepsilon\Vert x\Vert_{X_{M}}.
\]
Hence, according to Proposition \ref{Nonisomorphic}, $X$ and $X_{M}$ are not isomorphic.
\end{proof}

In the next two examples, we show that neither $\inf\varphi_{n}>0$ nor
$\sup\varphi_{n}<1$ is a necessary condition for $X$ and $X_{M}$ to be nonisomorphic.

\begin{example}
If $\theta<1$ and $\varphi_{n}=\frac{1}{n+1},$ then $X$ and $X_{M}$ are not isomorphic.
\end{example}

\begin{proof}
It suffices to show that $\left(  \theta_{n}\right)  $ satisfies $\left(
\lnot\dag\right)  ,$ $\left(  \ddag\right)  $ and $\lim_{N}\sup_{p}%
\frac{\Theta_{p}({N})}{\theta_{p}}=0.$ Note that
\[
\frac{\theta_{m+n}}{\theta_{n}}=\frac{n+1}{m+n+1}\theta^{m}.
\]
Hence
\[
\delta_{m}=\limsup_{n}\frac{\theta_{m+n}}{\theta_{n}}=\theta^{m} \to0
\]
as $m \to\infty$. Thus $\left(  \lnot\dag\right)  $ holds.

To see that $\left(  \theta_{n}\right)  $ satisfies $\left(  \ddag\right)  ,$
let $s_{1}<s_{2}<...<s_{R}$ be an arithmetic progression in $\mathbb{N}$. Note
that $s\mapsto\frac{s+1}{s+t+1}$ is a concave increasing function for
$s\geq0.$ Let $g\left(  s\right)  $ be the linear function interpolating
$\left(  s_{1},\frac{s_{1}+1}{s_{1}+t+1}\right)  $ and $\left(  s_{R}%
,\frac{s_{R}+1}{s_{R}+t+1}\right)  .$ Then
\begin{align*}
\sum_{i=1}^{R}\frac{\theta_{s_{i}+t}}{\theta_{s_{i}}}  &  =\theta^{t}%
\sum_{i=1}^{R}\frac{s_{i}+1}{s_{i}+t+1} \geq\theta^{t}\sum_{i=1}^{R}g\left(
s_{i}\right) \\
&  =\theta^{t}\frac{R}{2}\left[  g\left(  s_{1}\right)  +g\left(
s_{R}\right)  \right] \\
&  \quad\quad\text{since}\left(  g\left(  s_{i}\right)  \right)  _{i=1}%
^{R}\text{ is an arithmetic progression}\\
&  \geq\theta^{t}\frac{R}{2} \max\{g\left(  s_{1}\right)  , g\left(
s_{R}\right)  \} =\frac{R}{2}\max_{1\leq i\leq R}\frac{\theta_{s_{i}+t}%
}{\theta_{s_{i}}}.
\end{align*}
Hence $\left(  \ddag\right)  $ holds with $F\left(  R\right)  =\frac{2}{R}.$

Finally, if If $\left(  \ell_{i}\right)  _{i=1}^{N}$ is a sequence of positive
integers such that $\sum_{i=1}^{N}\ell_{i}=p,$ then at least one $\ell_{i}$ is
$\geq\frac{p}{N}.$ Without loss of generality, assume that $\ell_{1}\geq
\frac{p}{N}.$ Then
\[
\frac{1}{\ell_{1}+1} \leq\frac{N}{p+1}.
\]
Hence
\[
\prod_{i=1}^{N}\theta_{\ell_{i}} =\theta^{p}\prod_{i=1}^{N}\frac{1}{\ell
_{i}+1}\leq\theta^{p}( \frac{N}{p+1})( \frac{1}{2}) ^{N-1} =\frac{N}{2^{N-1}%
}\theta_{p}.
\]
Thus%
\[
\sup_{p}\frac{\Theta_{p}({N})}{\theta_{p}}\leq\frac{N}{2^{N-1}}.
\]
It follows from Proposition \ref{Nonisomorphic} and Theorem \ref{Th21} that
$X$ and $X_{M}$ are not isomorphic.
\end{proof}

\begin{example}
There exists a regular sequence $\left(  \theta_{n}\right)  $ with
$0<\theta<1$ and $\lim_{n}\varphi_{n}=1$ such that $X$ and $X_{M}$ are not isomorphic.
\end{example}

\begin{proof}
Let $0<\theta_{1}<\theta<1$ be given. Choose sequences $\left(  q_{n}\right)
$ and $\left(  K_{n}\right)  $ in $\mathbb{N}$ such that
\[
\theta^{q_{M+N+1}}\leq\frac{1}{24N^{2}}\theta_{1}^{2+s\left(  M,N\right)  }%
\]
and
\[
\frac{1}{K_{M+N+1}}\leq\frac{1}{144N^{2}}\theta_{1}^{3+s\left(  M,N\right)  }%
\]
if $0\leq M\leq N,$ where $s\left(  M,N\right)  =\sum_{i=1}^{M}K_{N+i}q_{N+i}$
if $0<M\leq N$ and $s\left(  0,N\right)  =0.$ Then choose a sequence $\left(
\varphi_{n}\right)  $ such that $\varphi_{1}=\frac{\theta_{1}}{\theta},$
$\left(  \varphi_{n}\right)  $ increases to $1$, $\varphi_{n+1}\leq
\frac{\varphi_{n}}{\theta}$ and $\lim_{N}\varphi_{s\left(  N,N\right)  }%
^{N}=0.$

Define $\theta_{n}=\varphi_{n}\theta^{n}.$ Then $\left(  \theta_{n}\right)  $
is a regular sequence such that $\lim\theta_{n}^{1/n}=\lim\varphi_{n}%
^{1/n}\theta=\theta.$ Since $\inf\varphi_{n}=\varphi_{1}>0,$ $\left(
\lnot\dag\right)  $ and $\left(  \ddag\right)  $ hold with $F\left(  R\right)
=\frac{1}{\varphi_{1}^{2}R}$ according to Proposition \ref{III}.

Given $N\in\mathbb{N}$, we claim that the sequences $\left(  p_{k}\right)
_{k=1}^{N}=\left(  q_{N+k}\right)  _{k=1}^{N}$ and $\left(  L_{k}\right)
_{k=1}^{N}=\left(  K_{N+k}\right)  _{k=1}^{N}$ satisfy conditions
\textbf{(A)}, \textbf{(B)}, and \textbf{(C)}. Indeed,
\[
\frac{\theta_{p_{M+1}+n}}{\theta_{n}}=\theta^{p_{M+1}}\frac{\varphi
_{p_{M+1}+n}}{\varphi_{n}}\leq\frac{\theta^{p_{M+1}}}{\varphi_{1}}\leq
\frac{\theta^{q_{N+M+1}}}{\theta_{1}}\leq\frac{1}{24N^{2}}\theta_{1}%
^{s(M,N)}.
\]
By regularity, $\theta_{n}\geq\theta_{1}^{n}.$ Hence
\[
\prod_{i=1}^{M}\theta_{L_{i}p_{i}}\geq\theta_{1}^{\sum_{i=1}^{M}L_{i}p_{i}%
}=\theta_{1}^{s(M,N)}%
\]
if $M>0$. Thus%
\[
\frac{\theta_{p_{M+1}+n}}{\theta_{n}}\leq\frac{1}{24N^{2}}\prod_{i=1}%
^{M}\theta_{L_{i}p_{i}}.
\]
Therefore, condition \textbf{(A)} is satisfied if $M>0$. If $M=0$, then
$s(M,N)=0$ and the vacuous product $\prod_{i=1}^{M}\theta_{L_{i}p_{i}}=1$ and
the result is clear.

To see that condition \textbf{(B)} is satisfied, we note that by the choice of
$(q_{n}),$ $q_{M+N+1}\geq2+s(M,N)$, which is equivalent to saying that
$p_{M+1}\geq2+\sum_{i=1}^{M}L_{i}p_{i}$ if $M>0$.

If $M>0$,
\begin{align*}
F(L_{M+1})  &  =\frac{1}{\varphi_{1}^{2}L_{M+1}}=\frac{1}{\varphi_{1}%
^{2}K_{M+N+1}}\leq\frac{1}{\theta_{1}^{2}\cdot144N^{2}}\theta_{1}^{3+s(M,N)}\\
&  \leq\frac{\theta_{1}}{144N^{2}}\theta_{1}^{\sum_{i=1}^{M}L_{i}p_{i}}%
\leq\frac{\theta_{1}}{144N^{2}}\prod_{i=1}^{M}\theta_{L_{i}p_{i}}.
\end{align*}
Therefore, condition \textbf{(C)} is also satisfied. Finally, we consider the
ratio%
\[
\frac{\Theta_{p\left(  r_{1},\dots,r_{N}\left(  k\right)  \right)  }}%
{\theta_{p\left(  r_{1},\dots,r_{N}\left(  k\right)  \right)  }}.
\]
If $\left(  \ell_{i}\right)  _{i=1}^{N}$ is a sequence in $\mathbb{N}$ such
that $\sum_{i=1}^{N}\ell_{i}=p\left(  r_{1},\dots,r_{N}\left(  k\right)
\right)  ,$ then
\begin{align*}
\prod_{i=1}^{N}\theta_{\ell_{i}}  &  =\theta^{p\left(  r_{1},\dots
,r_{N}\left(  k\right)  \right)  }\prod_{i=1}^{N}\varphi_{\ell_{i}}\\
&  \leq\theta^{p\left(  r_{1},\dots,r_{N}\left(  k\right)  \right)  }%
\varphi_{p\left(  r_{1},\dots,r_{N}\left(  k\right)  \right)  }^{N}%
\end{align*}
since $(\varphi_{n})$ is increasing and $0<\varphi_{n}<1$. Now
\begin{align*}
p\left(  r_{1},\dots,r_{N}\left(  k\right)  \right)   &  =r_{1}p_{1}%
+\cdots+r_{N-1}p_{N-1}+r_{N}\left(  k\right)  p_{N}\\
&  \leq L_{1}p_{1}+\cdots+L_{N}p_{N}=\sum_{i=1}^{N}K_{N+i}q_{N+i}=s\left(
N,N\right)  .
\end{align*}
Thus%
\begin{align*}
\prod_{i=1}^{N}\theta_{\ell_{i}}  &  \leq\theta^{p\left(  r_{1},\dots
,r_{N}\left(  k\right)  \right)  }\varphi_{s\left(  N,N\right)  }^{N}\\
&  =\frac{\varphi_{s\left(  N,N\right)  }^{N}}{\varphi_{p\left(  r_{1}%
,\dots,r_{N}\left(  k\right)  \right)  }}\theta_{p\left(  r_{1},\dots
,r_{N}\left(  k\right)  \right)  }\\
&  \leq\varphi_{s\left(  N,N\right)  }^{N}\varphi_{1}^{-1}\theta_{p\left(
r_{1},\dots,r_{N}\left(  k\right)  \right)  }.
\end{align*}
Hence
\[
\sup_{r_{1},\dots,r_{N-1}}\frac{\Theta_{p\left(  r_{1},\dots,r_{N}\left(
k\right)  \right)  }}{\theta_{p\left(  r_{1},\dots,r_{N}\left(  k\right)
\right)  }}\leq\varphi_{s\left(  N,N\right)  }^{N}\varphi_{1}^{-1}.
\]
since $\left(  \varphi_{N}\right)  $ is chosen such that $\lim_{N}%
\varphi_{s\left(  N,N\right)  }^{N}=0,$ we see that
\[
\lim_{N}\sup_{r_{1},\dots,r_{N-1}}\frac{\Theta_{p\left(  r_{1},\dots
,r_{N}\left(  k\right)  \right)  }}{\theta_{p\left(  r_{1},\dots,r_{N}\left(
k\right)  \right)  }}=0.
\]
Arguing as in the proof of Theorem \ref{Th16}, we may conclude that $X$ and
$X_{M}$ are not isomorphic.
\end{proof}

\end{document}